\documentclass[titlepage,11pt]{article}
\usepackage{latexsym}
\usepackage{amsfonts}
\usepackage{amsmath}
\usepackage{mathtools}
\usepackage{tikz}
\usepackage{comment}
\newcommand\blackslug{\hbox{\hskip 1pt \vrule width 4pt height 8pt depth 1.5pt
		\hskip 1pt}}
\newcommand\bbox{\hfill \quad \blackslug \bigbreak}
\def\DD{\hbox{-}}
\def\CC{\hbox{-}\cdots\hbox{-}}
\def\LL{,\ldots,}
\newcommand{\vare}{\varepsilon}

\newcommand{\dist}{\operatorname{dist}}

\newcommand{\cupcup}{\cup \cdots\cup}

\title{Asymptotic structure. II. Path-width and additive quasi-isometry}
\author{
	Tung Nguyen\thanks{Supported by a Porter Ogden Jacobus Fellowship, and AFOSR grant
		FA9550-22-1-0234, and by NSF grant  DMS-2154169.
		Current address: University of Oxford, Oxford, UK}\\
	Princeton University,\\ Princeton, NJ 08544, USA
	\and
	Alex Scott\thanks{Supported by EPSRC grant EP/X013642/1}\\
	University of Oxford, \\
	Oxford, UK
	\and
	Paul Seymour\thanks{Supported by AFOSR grant
		FA9550-22-1-0234, and by NSF grant DMS-2154169.}\\
	Princeton University,\\ Princeton, NJ 08544, USA}

\date{}
\date{}

\newtheorem{thm}{}[section]

\newcommand{\Proof}{\noindent{\bf Proof.}\ \ }

\begin{document}
	\maketitle
	\begin{abstract}
		We show that if a graph $G$ admits a quasi-isometry $\phi$ to a graph $H$ of bounded path-width, then 
		we can assign a non-negative integer length to each edge of $H$, such that the same function $\phi$ is a 
		quasi-isometry to this weighted version of $H$, with error only an additive constant. 
	\end{abstract}
	
	\section{Introduction}
	
	We need to begin with some definitions.
	Graphs in this paper may be infinite, and have no loops or parallel edges.
	If $X$ is a vertex of a graph $G$, or  
	a subset of the vertex set of $G$, or a subgraph of $G$, and the same for $Y$, then $\dist_G(X,Y)$ denotes the 
	distance in $G$
	between $X,Y$, that is, the number of edges in the shortest path of $G$ with one end in $X$ and the other in $Y$. (If no path exists we set $\dist_G(X,Y) = \infty$.)
	
	Let $G,H$ be graphs, and let $\phi:V(G)\to V(H)$ be a map.
	Let $L,C\ge 0$; we say that $\phi$ is an {\em $(L,C)$-quasi-isometry}, and {\em $G$ is $(L,C)$-quasi-isometric to $H$}, 
	if:
	\begin{itemize}
		\item for all $u,v$ in $V(G)$, if $\dist_G(u,v)$ is finite then $\dist_H(\phi(u),\phi(v))\le L \dist_G(u,v)+C$;
		\item for all $u,v$ in $V(G)$, if $\dist_H(\phi(u),\phi(v))$ is finite then $\dist_G(u,v)\le L \dist_H(\phi(u),\phi(v))+C$;
		and
		\item for every $y\in V(H)$ there exists $v\in V(G)$ such that $\dist_H(\phi(v), y)\le  C$.
	\end{itemize}
	What if we want $L=1$?
	There is a remarkable theorem of Chepoi, Dragan, Newman, Rabinovich, and Vaxès~\cite{chepoi}, also proved by Kerr~\cite{kerr}:
	\begin{thm}\label{kerr}
		For all $L,C$ there exists $C'$ such that if there is an $(L,C)$-quasi-isometry from a graph $G$ to a tree, then there
		is a $(1,C')$-quasi-isometry from $G$ to a tree.
	\end{thm}
	
	Is this special to trees, or can it be made much more general? 
	For instance, Agelos Georgakopoulos asked (in private communication) whether the same statement was true if we (twice) replace ``tree'' by ``planar graph''.
	Let $\mathcal{C}$ be a class of graphs. Under what 
	conditions on $\mathcal{C}$ can we say the following?
	
	``For all $L,C$ there exists $C'$ such that if there is an $(L,C)$-quasi-isometry from a graph $G$ to a member of $\mathcal{C}$, then there
	is a $(1,C')$-quasi-isometry from $G$ to a member of $\mathcal{C}$.''
	
	For this to be true, $\mathcal{C}$ must have some closure properties: for instance, if $H\in \mathcal{C}$ and $G$ is obtained from $H$
	by subdividing every edge once, there is a $(2,0)$-quasi-isometry from $G$ to $H$, but if we want there to be 
	a $(1,C')$-quasi-isometry from $G$ to a member of $\mathcal{C}$ then we need $\mathcal{C}$ to contain a graph much like $G$; and this is
	close to asking that $\mathcal{C}$ be closed under edge-subdivision.
	Similarly, if $H\in \mathcal{C}$ and $G$ is obtained from $H$ by contracting the edges in some matching of $H$, 
	there is a $(3,0)$-quasi-isometry from $G$ to $H$, and so we need $\mathcal{C}$ to be more-or-less closed under 
	edge-contraction. Is that enough, could the following be true?
	\begin{thm}\label{wildconj}
		{\bf Conjecture:} Let $\mathcal{C}$ be a class of connected graphs, closed under contracting edges and subdividing edges. 
		For all $L,C$ there exists $C'$ such that if there is an $(L,C)$-quasi-isometry from a graph $G$ to a member of $\mathcal{C}$, then there
		is a $(1,C')$-quasi-isometry from $G$ to a member of $\mathcal{C}$.
	\end{thm}
	For instance, if $G,H$ are respectively the infinite square lattice and the infinite triangular lattice, there is a quasi-isometry between them,
	but no $(1,C)$-quasi-isometry (for any constant $C$); but there is a $(1,2)$-quasi-isometry from $G$ to a graph obtained by subdividing edges
	of $H$, and a $(1,100)$-quasi-isometry from $H$ to a graph obtained by subdividing and contracting edges
	of $G$ (we omit the proofs of all these statements).
	
	We are far from proving the conjecture \ref{wildconj} in general, but we will prove a special case, which we will explain next.
	A {\em path-decomposition} of a graph $G$ is (possibly infinite or doubly-infinite) sequence $(B_t:t\in T)$, where
	$T$ is a set of integers, and $B_t$
	is a subset of $V(G)$ for each $t\in V(T)$ (called a {\em bag}), such that:
	\begin{itemize}
		\item $V(G)$ is the union of the sets $B_t\;(t\in T)$;
		\item for every edge $e=uv$ of $G$, there exists $t\in T$ with $u,v\in B_t$; and
		\item for all $t_1,t_2,t_3\in T$, if $t_1\le t_1\le t_3$, then
		$B_{t_1}\cap B_{t_3}\subseteq B_{t_2}$.
	\end{itemize}
	The {\em width} of a path-decomposition $(T,(B_t:t\in V(T)))$ is the maximum of the numbers $|B_t|-1$ for $t\in V(T)$,
	or $\infty$ if there is no finite maximum;
	and the {\em path-width} of $G$ is the minimum width of a path-decomposition of $G$.
	
	We will prove: 
	\begin{thm}\label{pathwidth}
		For all $L,C,k$ there exists $C'$ such that if there is an $(L,C)$-quasi-isometry from a graph $G$ to a graph $H$ with path-width at most $k$, then there
		is a $(1,C')$-quasi-isometry from $G$ to a graph $H'$ obtained from $H$ by subdividing and contracting edges.
	\end{thm}
	In fact $C'$ can be taken to be $\max(L,C)^{O(k)}$.

	Let $\mathbb{N}$ denote the set of nonnegative integers.
	Let $H$ be a graph and let $w:E(H)\to\mathbb{N}$ be some function; we call $(H,w)$  a {\em weighted graph}. One can define quasi-isometry for
	weighted graphs in the natural way, defining $\dist_{(H,w)}(u,v)$ to be the minimum of $w(P)$ over all paths $P$ of $H$ between $u,v$,
	where $w(P)$ means $\sum_{e\in E(P)}w(e)$. Subdividing and contracting edges of $H$ is closely related to moving from $H$ to $(H,w)$ 
	for an appropriate $w$, so we could express \ref{pathwidth} in terms of weighted graphs. In this modified form of \ref{pathwidth},
	rather than replacing $H$ by $H'$, we keep $H$ and just put weights on its edges.
	But something much stronger is true: we don't need to change the quasi-isometry either.
	
	\begin{thm}\label{weightedpathwidth}
		For all $L,C,k$ there exists $C'$ such that if $\phi$ is an $(L,C)$-quasi-isometry from a graph $G$ to a graph $H$ with
		path-width at most $k$, then there is a function $w:E(H)\to \mathbb{N}$ such that the same function $\phi$
		is a $(1,C')$-quasi-isometry from $G$ to  the weighted graph $(H,w)$.
	\end{thm}
	
	In view of the conjecture \ref{wildconj}, one might ask whether the path-width restriction is necessary.
	It cannot just be omitted, because Davies, Hatzel and Hickingbotham~\cite{davies} very recently showed the following:
	\begin{thm}\label{daviesthm}
		For every integer $C>0$, and and every real $\vare>0$, there exist graphs $G$ and $H$ such that $G$ is
		$(1 + \vare, 1)$-quasi-isometric to $H$ but there is no map $w:E(H) \to \mathbb{R}$ of $H$ such
		that $G$ is $(1, C)$-quasi-isometric to $(H, w)$.
	\end{thm}
	(Curiously, this does not disprove the conjecture \ref{wildconj}.)
	It still might be true that we can replace ``with path-width at most $k$'' in \ref{weightedpathwidth} by something less restrictive,
	for instance, by ``with tree-width at most $k$'', or ``that is planar'', or indeed ``that is in $\mathcal{C}$''
	where $\mathcal{C}$ is any minor-closed class of graphs that does not contain all finite graphs.
	
	Is \ref{wildconj} true at least when $\mathcal{C}$ is the class of connected graphs with tree-width at most $k$? (A closely-related question
	was considered by Dragan and Abu-Ata~\cite{dragan}.) Yes when $k=1$, by
	\ref{kerr}, and indeed one can show that \ref{pathwidth} also holds in this case (see the proof of \ref{kerr} in~\cite{shortcuts}).
	What about tree-width two? A special case is when $\mathcal{C}$ is the class of all connected outerplanar graphs, and we can prove
	\ref{wildconj} in that case. (A hint for the proof: every connected outerplanar graph is quasi-isometric to a graph in which every non-trivial block is a cycle.) But for
	tree-width two in general, the result is open, as is the following weaker statement:
	
	\begin{thm}\label{tw2conj}
		{\bf Conjecture: }For all $L,C$ there exist $C',k$ such that if there is an $(L,C)$-quasi-isometry from a graph $G$ to a graph of tree-width at most two,
		then there is a $(1,C')$-quasi-isometry from $G$ to a graph of tree-width at most $k$.
	\end{thm}
	
	So far, our statements are true for infinite graphs as well as for finite graphs, but we want to make an adjustment,
	because  path-width is not the ``right'' concept for infinite graphs. A graph has tree-width at most $k$ if and only of all its finite
	subgraphs have tree-width at most $k$ (Thomas~\cite{robin}), but the same is not true for path-width. For instance, the graph consisting of the disjoint union of infinitely many one-way infinite paths has infinite path-width, and so does the disjoint union of infinitely many 
	copies of the infinite ``star'' (one vertex with countably many neighbours); and so does any graph with uncountably many vertices and 
	no edges. There is a more appropriate concept. Let us say a {\em line} is a set that is linearly ordered by some relation $<$; and 
	a {\em line-decomposition} is a family $(B_t:t\in T)$, where $T$ is a line, satisfying the same three conditions as in the definition 
	of path-decomposition. We define the {\em width} of a line-decomposition to be the maximum of $|B_t|-1$  over all $t\in T$ if this exists,
	and otherwise the width is infinite. The {\em line-width} of $G$ is the minimum integer $k$ such that $G$ admits a line-decomposition 
	of width at most $k$, if this exists, and otherwise the line-width is infinite. For finite graphs, path-width and line-width are the same, but 
	for infinite graphs, they may be different (for instance, in the three examples above), and line-width behaves better. We proved in~\cite{subtrees}
	that a graph has line-width at most $k$ if and only if all its finite subgraphs have path-width at most $k$.
	
	All the theorems about path-width mentioned so far are also true for line-width, and expressing them this way makes them stronger and more general.
	In particular, we will prove:
	\begin{thm}\label{weightedlinewidth}
		For all $L,C,k$ there exists $C'$ such that if $\phi$ is an $(L,C)$-quasi-isometry from a graph $G$ to a graph $H$ with
		line-width at most $k$, then there is a function $w:E(H)\to \mathbb{N}$ such that the same function $\phi$
		is a $(1,C')$-quasi-isometry from $G$ to  the weighted graph $(H,w)$.
	\end{thm}

	Here is an application. A. Georgakopoulos in private communication showed that  for all $L,C$ there exists $C'$ such that
	if a finite graph $G$ is $(L,C)$-quasi-isometric to a cycle, then $G$ is $(1,C')$-quasi-isometric to a cycle. 
	This immediately follows from \ref{pathwidth}. Similarly, we (unpublished) proved some time ago the following result about fat minors 
	(we omit the definitions of fat minor, since we will not need them any more in this paper): for all $k,C$, there exists 
	$C'$ such that if a graph $G$ does not contain $K_{1,k}$ as a $C$-fat minor, then there is a $(1,C')$-quasi-isometry from $G$ to a graph 
	not containing $K_{1,k}$ as a minor. This strengthened a result of Georgakopoulos and Papasoglu~\cite{agelos} that all $k,C$, there exist
	$L, C'$ such that if $G$ does not contain $K_{1,k}$ as a $C$-fat minor, then there is an $(L,C')$-quasi-isometry from $G$ to a graph
	not containing $K_{1,k}$ as a minor.
	Our proof was complicated, but connected graphs with no $K_{1,k}$ minor are have line-width at most $k-1$,
	and so our result follows via \ref{weightedlinewidth} from that 
	of Georgakopoulos and Papasoglu.


	\section{Finding a weighting in the neighbourhood of $\phi(P)$}

	If $(H,w)$ is a weighted graph,
	the {\em size} of $w$ is  the maximum of $w(e)$ over all $e\in E(G)$, assuming this exists: we will only use
	weighted graphs with bounded size.
	
	Let us reiterate a definition, more explicitly. Let $G$ be a graph and let $(H,w)$ be a weighted graph.
	A map $\phi$ from $V(G)$ to $V(H)$ is an {\em $(L,C)$-quasi-isometry from $G$ to $(H,w)$}
	if:
	\begin{itemize}
		\item for all $u,v$ in $V(G)$, if $\dist_G(u,v)$ is finite then $\dist_{(H,w)}(\phi(u),\phi(v))\le L \dist_G(u,v)+C$;
		\item for all $u,v$ in $V(G)$, if $\dist_{(H,w)}(\phi(u),\phi(v))$ is finite then $\dist_G(u,v)\le L \dist_{(H,w)}(\phi(u),\phi(v))+C$; and
		\item for every $y\in V(H)$ there exists $v\in V(G)$ such that $\dist_{(H,w)}(\phi(v), y)\le  C$.
	\end{itemize}
	
	For inductive purposes, it is more convenient to prove the following stronger version of \ref{weightedlinewidth}:
	\begin{thm}\label{weightedline}
		Let $L,C,k\ge 0$ be integers; then there exist $C', W$ with the following property. Let $H$ be a graph with line-width at most $k$,
		and let $\phi$ be an $(L,C)$-quasi-isometry from a graph $G$ to $H$. Then 
		there is a function $w: E(H)\to \mathbb{N}$ with size at most $W$ such that $\phi$ is a $(1,C')$-quasi-isometry from $G$ to $(H,w)$.
	\end{thm}
	
	Instead of working with $(L,C)$-quasi-isometries, we could replace $L,C$ by their common maximum, and so it would be enough to work with $(C,C)$-quasi-isometries. Actually, we prefer to use $(C-1,C)$-quasi-isometries, because then the small terms in the various numerical
	expressions that come up are easier to dispose of.
	
	If $P$ is a path and $u,v\in V(P)$, we denote by $P[u,v]$ the subpath between $u,v$.
	A {\em geodesic} in a graph $G$ means a path $P$ of $G$ (possibly infinite) such that for every two vertices $u,v\in V(P)$, the 
	subpath $P[u,v]$
	is a shortest path of $G$ between $u,v$.
	If $(H,w)$ is a weighted graph, a 
	{\em $w$-geodesic} of $H$ means a path $P$ of $H$ such that $\dist_{(H,w)}(u,v) = w(P[u,v])$ for all $u,v\in V(P)$.
	An {\em integer interval} means a set of integers $I$, finite or infinite, such that if $i,k\in I$ and $j$ is an integer with $i<j<k$
	then $j\in I$.
	
	Let us sketch an outline of the proof of \ref{weightedline}. We work by induction on the line-width. Let $(B_t:t\in T)$
	be a line-decomposition of $H$ of width at most $k$. Thus, $H$ can be thought of intuitively as a long, thin graph in some sense, 
	and so is $G$. One would expect there to be a geodesic $P$ of $G$ running through all the preimages of the bags $B_t$. 
	If we have such a geodesic, then the vertices $\phi(p)$ move along the length of $H$ as $p$ runs through the vertices of $P$;
	and we find a path $Q$ of $H$ staying close to all these vertices. Since $P$ is a geodesic in $G$, we can arrange weights $w_1$
	in $H$ such that $Q$ is a $w_1$-geodesic in $H$, and in addition, for every two vertices of $P$, their distance in $G$ is about the same
	(up to an additive error) as the weighted distance between their images in $H$. 
	One can show that, then, the same is true for all vertices $u,v$ of $G$ that are at most a constant distance 
	from $P$: $\dist_G(u,v)$ and $\dist_{(H,w_1)}(\phi(u), \phi(v))$ differ only by a constant.
	We still need to work 
	on the set $B$ of vertices that are far from $P$; 
	but they map to a subgraph of $H$
	with line-width less than $k$, so we can use the inductive hypothesis for them, provided that the restriction of $\phi$ to $B$
	is still a quasi-isometry with bounded parameters. To fix this last condition needs some fiddling around; we have to add some
	new vertices to $B$ to make the distances in $B$ the same as they were in $G$. But this works.
	
	There is a major difficulty in finding the geodesic $P$. It is easy to obtain if $H$ is finite, but if $H$ is infinite it needs 
	a lot more work. (One difficulty is that such a geodesic need not exist: we have to grow $G$ into a bigger graph to obtain $P$.)
	We have arranged the paper with the arguments to obtain $P$ at the end, so the reader who wants to understand the proof for 
	finite $H$ does not have to wade through pages of argument for infinite graphs.
	
	If $I,J$ are nonempty sets of integers, we say that $J$ is {\em cofinal} with $I$ if either they both have a maximum element and these
	elements are equal, or neither has a maximum element; and either they both have a minimum element and these elements are equal, or
	neither has a minimum element.
	Our objective in this section is, given the geodesic $P$ in $G$, to find an appropriate path $Q$ of $H$ as above, and find the 
	weight function $w_1$ that makes $Q$ a geodesic and spaces appropriately the images of the vertices in $P$ . More exactly:
	\begin{thm}\label{fixgeo0}
		Let $C\ge 4$, and let $\phi$ be a $(C-1,C)$-quasi-isometry from a graph
		$G$ to a graph $H$.
		Let $P$ be a geodesic in $G$, with vertices $p_i\;(i\in I)$ in order, where $I$ is an integer interval.
		Then there is a function $w:E(H)\to \mathbb{N}$, with size at most $32C^4$, and a path $Q$ of $H$, and $J\subseteq I$ cofinal with $I$,
		and distinct vertices 
		$r_j\;(j\in J)$, in order in $Q$, with the following properties:
		\begin{itemize}
			\item $Q$ is a $w$-geodesic in $H$;
			\item $\dist_{(H,w)}(r_i,r_j) = j-i$ for all $i,j\in J$ with $i<j$;
			\item for all $i\in I$ there exists $j\in J$ with $|j-i|\le C^2$ and $\dist_H(\phi(p_i), r_j)< C^3$; 
			\item $\dist_H(\phi(p_j), r_j)\le 2C$ for each $j\in J$; and
			\item for each $v\in V(Q)$ there exists $j\in J$ such that $\dist_H(v, \phi(p_j))\le C$.
		\end{itemize}
	\end{thm}
	
	We divide the proof into two steps. First, we show:
	\begin{thm}\label{getweights}
		Let $L\ge 1$ be an integer. Let $J$ be a set of integers, 
		let $Q$ be a path of a graph $H$, and let $r_i\in V(Q)$ for each $i\in J$, all distinct and numbered in order on $Q$.
		Suppose that:
		\begin{itemize}
			\item $Q$ is the union of the subpaths $Q[r_i,r_j]$ for $i,j\in J$;
			\item for all $i,j\in J$ with $i\le j$, if none of $i+1\LL j-1$ belong to $J$, then $j-i\le L$;
			\item $\dist_H(r_i,r_j)\ge (j-i)/L$ for all $i,j\in J$ with $j>i$; and
			\item $\dist_Q(r_i,r_j)\le L(j-i)$ for all $i,j\in J$ with $j>i$.
		\end{itemize}
		Then there is a function $w:E(H)\to \mathbb{N}$ with size at most $L(2L+1)$, such that 
		$Q$ is a $w$-geodesic of $H$, and 
		$\dist_{(H,w)}(r_i,r_j)=j-i$ for all $i,j\in J$ with $j\ge i$.
	\end{thm}
	\Proof
	We may assume that $|J|>1$. Let $R=\{r_i:i\in J\}$.
	Let us say a {\em gap} is a subpath of $Q$ of length at least one, with both ends in $R$ and with no internal vertex in this 
	set. Thus all gaps have length at most $L$, and every vertex of $Q$ belongs to a gap. By hypothesis, the vertices 
	$r_i\;(i\in J)$ are numbered in their order in $Q$. This extends to an ordering of the vertex set of $Q$, which we call ``later than''.
	We begin with:
	\\
	\\
	(1) {\em If $x,y\in V(Q)$. Then $Q[x,y]$ is contained in $Q[r_i,r_j]$ for some $i,j\in J$ with 
		$0\le j-i\le L(2L+1)\dist_H(x,y)$.}
	\\
	\\
	We may assume that $y$ is later than $x$. Choose $i\in J$ maximum such that $x$ is later than or equal to $r_i$, and choose $j\in J$
	minimum such that $r_j$ is later than or equal to $y$. 
	Thus $i<j$, and  
	$\dist_H(x,y)\ge \dist_H(r_i,r_j)-2L$. But $\dist_H(r_i,r_j)\ge (j-i)/L$,
	and so $\dist_H(x,y)\ge (j-i)/L-2L$. Since $x\ne y$, it follows that $\dist_H(x,y)\ge 1$ and so 
	$2L\dist_H(x,y)\ge 2L$. Consequently, $(2L+1)\dist_H(x,y) \ge (j-i)/L$. This proves (1).
	
	\bigskip
	
	For each gap $Q[r_i,r_j]$, and each edge $e$ of this gap, define $w(e) = j-i$ if $e$ is incident with $r_i$, and $w(e) = 0$ otherwise. 
	It follows that $w(Q[r_i,r_j]) = j-i$ for all $i,j\in J$ with $i\le j$.
	Define $w(e) = L(2L+1)$ for every edge $e$ of $H$ not in $E(Q)$. 
	
	It remains to show that 
	$Q$ is a $w$-geodesic of $H$.
	To show this, let $x,y\in V(Q)$, and let $P$
	be a $w$-geodesic of $H$ between $x,y$. We need to show that $w(P)\ge w(Q[x,y])$, and we prove this by induction on $|E(P)|$. If some internal vertex $z$ of $P$ belongs to $V(Q)$, then from the inductive hypothesis, $w(P[x,z])\ge w(Q[x,z])$, and $w(P[z,y])\ge w(Q[z,y])$,
	and adding, it follows that 
	$$w(P)\ge  w(Q[x,z])+w(Q[z,y])\ge w(Q[x,y])$$ 
	as required. So we may assume that no internal vertex of $P$
	belongs to $V(Q)$. We may assume that $P\ne Q[x,y]$, and so no edge of $P$ is in $E(Q)$, and therefore 
	$$w(P) = L(2L+1)|E(P)|\ge L(2L+1)\dist_H(x,y).$$
	By (1), $Q[x,y]$ is contained in $Q[r_i,r_j]$ for some $i,j\in J$ with 
	$0\le j-i\le L(2L+1)\dist_H(x,y).$
	Since $w(Q[x,y])\le w(Q[r_i,r_j])= j-i$, we deduce that 
	$$w(Q[x,y])\le L(2L+1)\dist_H(x,y)\le w(P).$$ 
	This proves \ref{getweights}.~\bbox
	
	The second step is:
	\begin{thm}\label{mapgeo}
		Let $C\ge 2$, and let $\phi$ be a $(C-1,C)$-quasi-isometry from a graph
		$G$ to a graph $H$.
		Let $P$ be a geodesic in $G$, with vertices $p_i\;(i\in I)$, numbered in order, where $I$ is an integer interval.
		Then there exist $J\subseteq I$, cofinal with $I$, and a path $Q$ of $H$, and distinct vertices $r_j\;(j\in J)$ of $Q$, 
		in order in $Q$,
		with the following properties:
		\begin{itemize}
			\item $Q$ is the union of the subpaths $Q[r_i,r_j]$ for $i,j\in J$;
			\item for all $i,j\in J$ with $i\le j$, if none of $i+1\LL j-1$ belong to $J$, then $j-i\le 2C^2$;
			\item $\dist_H(r_i,r_j)\ge (j-i)/(4C^2-1)$ for all $i,j\in J$ with $j>i$; 
			\item $\dist_Q(r_i,r_j)\le 2C(j-i)$ for all $i,j\in J$ with $j>i$; 
			\item $\dist_H(\phi(p_j), r_j)\le 2C$ for each $j\in J$; and
			\item for each $v\in V(Q)$, there exists $j\in J$ such that $\dist_H(v,\phi(p_j))\le C$.
		\end{itemize}
	\end{thm}
	\Proof
	We may assume that $|I|\ge 1$. There are four cases, depending whether $I$ is finite, or one-way infinite (in two possible ways) or two-way infinite. 
	Suppose first that $I$ has a minimum; then by renumbering, we can assume this minimum is zero. Let $i_0=0$. Inductively, 
	having defined $i_0\LL i_k\in I$, with $i_0<\cdots < i_k$, if $i_k$ is the maximum of $I$, stop. Otherwise let $i_{k+1}\in I$
	be maximum such that $\dist_H(\phi(p_{i_k}), \phi(p_{i_{k+1}}))\le 2C$; and let $T_{k}$ be a geodesic between 
	$\phi(p_{i_k}),\phi(p_{i_{k+1}})$. This exists and $i_{k+1}>i_k$, since 
	$\dist_H(\phi(p_{i_k}), \phi(p_{1+i_{k}}))\le 2C$ (because $\phi$ is a $(C-1,C)$-quasi-isometry). 
	Thus $T_0,T_1\LL $ are all paths in $H$ of length at most $2C$. Certainly $T_k$ meets $T_{k+1}$ for each $k$, since they share an 
	end-vertex and perhaps more. We claim that $T_h, T_k$ are vertex-disjoint if $k\ge h+2$; because suppose $v$ is a vertex in both paths.
	The sum of the distances between $v$ and $\phi(p_{i_{h-1}}), \phi(p_{i_{h}}), \phi(p_{i_{k-1}}),\phi(p_{i_{k}})$  is at most $4C$,
	since the sum of the first two is the length of $T_h$, and the last two sum to $T_k$. Consequently, either there is a path
	between $\phi(p_{i_{h-1}}),\phi(p_{i_{k-1}})$ of length at most $2C$, or one between $\phi(p_{i_{h}}),\phi(p_{i_{k}})$ of length at most $2C$,
	and this contradicts the definition of $i_{h}$ or of $i_{h+1}$.  Thus, in the sequence $T_0,T_1,\ldots$, non-consecutive terms are 
	vertex-disjoint. Let $Q$ be the path defined as follows: start with the subpath of $T_0$ from $\phi(p_0)$ to the first vertex of $T_0$
	in $T_1$; then follow $T_1$ to the first vertex of $T_1$ in $T_2$; and so on, for each integer $i\ge 0$ if $I$ is infinite,
	or until $i$ is the maximum element of $I$. In the second case, let $i_k$ be this maximum element: extend $Q$ along $T_{k-1}$ to 
	$\phi(p_{i_k})$.

	Let 
	$J=\{i_0,i_1,\ldots \}$. Let $r_0 = \phi(p_0)$, and if $I$ has a maximum element $i_k$ let $r_{i_k} = \phi(p_{i_k})$. For each 
	$i_k\in J$ not the minimum or maximum element of $I$, choose $r_{i_k}\in V(Q)\cap V(T_{k-1})\cap V(T_k)$. 
	It follows that the vertices $r_j (j\in J)$ are distinct and in order in $Q$,
	and $Q$ is the union of the subpaths $Q[r_i,r_j]$ for $i,j\in J$. We claim:
	\\
	\\
	(1) {\em The following hold:
		\begin{itemize}
			\item for all $i,j\in J$ with $i\le j$, if none of $i+1\LL j-1$ belong to $J$, then $j-i\le 2C^2-C$;
			\item $\dist_H(r_i,r_j)\ge (j-i)/(4C^2-1)$ for all $i,j\in J$ with $j>i$; 
			\item $\dist_Q(r_i,r_j)\le 2C(j-i)$ for all $i,j\in J$ with $j>i$;
			\item for each $j\in J$, $\dist_H(\phi(p_j), r_j)\le 2C$; and 
			\item for each $v\in V(Q)$ there exists $j\in J$ with $\dist_H(v,\phi(p_j))\le C$. 
		\end{itemize}
	}
	\noindent
	For the first bullet, let $i,j\in J$ with $i\le j$, such that none of $i+1\LL j-1$ belong to $J$. We may assume that 
	$j>i$ and it follows that $i=i_k$ and $j=i_{k+1}$ for some choice of $k$. Hence $T_k$ exists, and joins $\phi(p_{i_k}),\phi(p_{i_{k+1}})$.
	Consequently $\dist_H(\phi(p_{i_k}),\phi(p_{i_{k+1}}))\le 2C$, and so 
	$$\dist_G(p_{i_k},p_{i_{k+1}})\le (C-1)(2C)+C= 2C^2-C.$$
	
	For the second bullet, let $i= i_k$ and $j= i_\ell$ where $\ell>k$. Then 
	$$\dist_H(r_i, r_j)\ge \dist_H(\phi(p_{i_k}), \phi(p_{i_\ell}))-4C.$$
	Since 
	$$j-i=i_\ell-i_k=\dist_G(p_{i_k},p_{i_\ell})\le (C-1) \dist_H(\phi(p_{i_k}), \phi(p_{i_\ell}))+C,$$
	we deduce that 
	$$\dist_H(r_i, r_j)\ge (j-i-C)/(C-1)-4C=(j-i)/(C-1) -(C/(C-1)+4C).$$
	Since also $\dist_H(r_i, r_j)\ge 1$, it follows that 
	$$\dist_H(r_i, r_j)\ge (j-i)/(C-1)-(C/(C-1)+4C)\dist_H(r_i, r_j),$$ 
	and so 
	$$\dist_H(r_i, r_j)\ge \frac{j-i}{(C-1)(1+C/(C-1)+4C)}\ge (j-i)/(4C^2-1),$$
	as claimed.
	
	For the third bullet, again let $i= i_k$ and $j= i_\ell$ where $\ell>k$. The subpath $Q[r_{i_k},r_{i_\ell}]$ is the union of subpaths
	of $T_k, T_{k+1}\LL T_{\ell-1}$, and so has length at most $2C(\ell-k)\le 2C(i_{\ell}-i_k)$, since $\ell-k\le i_{\ell}-i_k$.
	This proves the third bullet.
	
	For the fourth bullet, let $j=i_k\in J$; then $r_j, \phi(p_j)$ are both vertices of $T_k$, which has length at most $2C$. 
	Finally, for the fifth bullet, let $v\in V(T_k)$; then $T_k$ has ends $\phi(p_{i_k}),\phi(p_{i_{k+1}})$, and so $v$ has distance in $H$ at most $C$ from one of these ends.
	This proves (1).
	
	\bigskip
	So in the case when $I$ has a minimum, the theorem holds. 
	Thus we may assume that $I$ has no minimum, and similarly that it has no maximum; and so $I=\mathbb{Z}$. 
	Choose integers $i_0<i_1$ with the interval $[i_0,i_1]$ maximal such that $\dist_H(\phi(p_{i_0}),\phi(p_{i_1}))\le 2C$, 
	and let $T_0$ be a geodesic between $\phi(p_{i_0}),\phi(p_{i_1})$.
	We define $i_1, i_2\LL$ inductively as before; that is, for each $k\ge 1$, let $i_{k+1}\in I$
	be maximum such that $\dist_H(\phi(p_{i_k}), \phi(p_{i_{k+1}}))\le 2C$, and let $T_{k}$ be a geodesic between  
	$\phi(p_{i_k}),\phi(p_{i_{k+1}}))$. Define the 1-way infinite path (previously called $Q$) as before, and let us call it $Q^+$.
	
	Now we define $i_{-1}, i_{-2}$ and so on, inductively:
	having defined $i_0, i_{-1}\LL i_{-k}$, let $i_{-k-1}$ be minimum such that 
	$\dist_H(\phi(p_{i_{-k}}), \phi(p_{i_{-k-1}}))\le 2C$; and let $T_{-k}$ be a geodesic between  
	$\phi(p_{i_{-k}}),\phi(p_{i_{-k-1}})$. Let $Q^{-}$ be defined in the same way that we defined $Q^+$. Now if $h<0$ and $j>0$,
	the paths $T_h, T_j$ are vertex-disjoint: because if they meet, then as before, 
	either there is a path
	between $\phi(p_{i_{h-1}}),\phi(p_{i_{k-1}})$ of length at most $2C$, or one between $\phi(p_{i_{h}}),\phi(p_{i_{k}})$ of length at most $2C$, in either case contrary to the maximality of the interval $[i_0,i_1]$. So $Q^+,Q^-$ meet only in vertices of $T_0$;
	and hence there is a path $Q$ contained in $Q^+\cup Q^-$, including all of $Q^+$ and all of $Q^-$ except possibly for some vertices in
	$T_0$, and containing at least one vertex of $T_0$. Then as before $Q$ satisfies the theorem. This proves \ref{mapgeo}.~\bbox
	
	By combining these two results, we obtain \ref{fixgeo0}, which we restate:
	\begin{thm}\label{fixgeo}
		Let $C\ge 2$, and let $\phi$ be a $(C-1,C)$-quasi-isometry from a graph
		$G$ to a graph $H$.
		Let $P$ be a geodesic in $G$, with vertices $(p_i\;(i\in I)$ in order, where $I$ is an integer interval.
		Then there is a function $w:E(H)\to \mathbb{N}$, with size at most $32C^4$, and a path $Q$ of $H$, and $J\subseteq I$, and distinct vertices
		$r_j\;(j\in J)$, in order in $Q$, with the following properties:
		\begin{itemize}
			\item $Q$ is a $w$-geodesic in $H$;
			\item $\dist_{(H,w)}(r_i,r_j) = j-i$ for all $i,j\in J$ with $i<j$;
			\item for all $i\in I$ there exists $j\in J$ with $|j-i|\le C^2$ and $\dist_H(\phi(p_i), r_j)< C^3$;
			\item $\dist_H(\phi(p_j), r_j)\le 2C$ for each $j\in J$; and
			\item for each $v\in V(Q)$ there exists $j\in J$ such that $\dist_H(v, \phi(p_j))\le C$.
		\end{itemize}
	\end{thm}
	\Proof
	By \ref{mapgeo}, 
	there exist $J\subseteq I$, and a path $Q$ of $H$, and distinct vertices $r_j\;(j\in J)$ of $Q$,
	in order in $Q$,
	with the following properties:
	\begin{itemize}
		\item $J$ is cofinal with $I$, and $Q$ is the union of the subpaths $Q[r_i,r_j]$ for $i,j\in J$;
		\item for all $i,j\in J$ with $i\le j$, if none of $i+1\LL j-1$ belong to $J$, then $j-i\le 2C^2$;
		\item $\dist_H(r_i,r_j)\ge (j-i)/(4C^2-1)$ for all $i,j\in J$ with $j>i$;
		\item $\dist_Q(r_i,r_j)\le 2C(j-i)$ for all $i,j\in J$ with $j>i$;
		\item $\dist_H(\phi(p_j), r_j)\le 2C$ for each $j\in J$; and
		\item for each $v\in V(Q)$, there exists $j\in J$ such that $\dist_H(v,\phi(p_j))\le C$.
	\end{itemize}
	By \ref{getweights}, taking $L=4C^2-1$,
	there is a function $w:E(H)\to \mathbb{N}$ with size at most $L(2L+1)\le 32C^4$, such that
	$Q$ is a $w$-geodesic of $H$, and
	$\dist_{(H,w)}(r_i,r_j)=j-i$ for all $i,j\in J$ with $j\ge i$.
	Thus the first, second, fourth and fifth bullets of the theorem hold.
	For the third, 
	let $i\in I$. There exist $j_1,j_2\in J$ with $j_1\le i\le j_2$ such that $j_2-j_1\le 2C^2$, since $J$ is cofinal with $I$, 
	and from the second bullet above. So there exists $j\in J$
	with $\dist_G(p_i,p_j) = |j-i|\le C^2$. Consequently $\dist_{H}(\phi(p_i),\phi(p_j))\le (C-1)C^2+C$.
	Since $\dist_H(\phi(p_j), r_j)\le 2C$ from the fifth bullet above, it follows that 
	$\dist_{H}(\phi(p_i),r_j)\le (C-1)C^2+3C<C^3$. This proves the third bullet is satisfied, and so proves \ref{fixgeo}.~\bbox

	
	\section{Extending the local weighting to the whole of $H$}
	Now we turn to the second part of the proof of \ref{weightedline}.
	We have found the geodesic $P$ of $G$, and the path $Q$ of $H$ and a weight function $w_1$ that makes $Q$ a $w_1$-geodesic and makes $\phi$
	have only additive error for vertices close to $P$; and we know that the theorem is true for graphs $H$ of smaller line-width. 
	We want to redefine the weights on edges of $H$ far from $Q$, to obtain a weight function $w$ on $H$ that satisfies \ref{weightedline}.
	
	Let us say a function $\kappa:\mathbb{N}\to \mathbb{N}$ is an {\em additive bounder} for a class $\mathcal{C}$ of graphs if
	for all $C\ge 1$, and every $(C-1,C)$-quasi-isometry $\phi$ from a graph $G$ to a graph 
	$H\in \mathcal{C}$,
	there is a function $w:E(H)\to \mathbb{N}$ with size at most $\kappa(C)$, such that $\phi$
	is a $(1,\kappa(C))$-quasi-isometry from $G$ to $(H,w)$. 
	
	A class $\mathcal{C}$ of graphs is {\em hereditary} if for every $H\in \mathcal{C}$, all induced subgraphs of $H$ also belong to 
	$\mathcal{C}$. The next result is the second step of the proof of \ref{weightedline}. (The additive bounder and the hereditary class in the
	statement are just a way to avoid talking about the induction on line-width. When we apply this result, $\mathcal{C}$ will be the 
	class of all graphs with line-width at most $k-1$, and $\kappa(C)$ will be a value of $C'$ that satisfies \ref{weightedline} with 
	$L=C-1$ and with $k$ replaced by $k-1$.) 
	\begin{thm}\label{usegeo}
		Let $\mathcal{C}$ be a hereditary class of graphs, with an additive bounder $\kappa$. For all $c\ge 2$ there exists $c_0$ with the following property. Suppose that:
		\begin{itemize}
			\item $\phi$ is a $(c-1,c)$-quasi-isometry from a graph $G$ to a graph $H$; 
			\item 
			$P$ is a geodesic in $G$, with vertices $p_i\;(i\in I)$ in order, where $I$ is an interval of integers, and $Q$ is a path of $H$;
			\item $J\subseteq I$, cofinal with $I$, and $r_j\;(j\in J)$ are vertices of $Q$ in order, 
			and $Q$ is the union of the subpaths $Q[r_i,r_j]$ for $i,j\in J$;
			\item $\dist_H(\phi(p_j), r_j)\le c$ for each $j\in J$, and for all $i\in I$ there exists $j\in J$ with $|j-i|\le c$ and $\dist_H(\phi(p_i), r_j)< c$;
			\item for each $v\in V(Q)$ there exists $j\in J$ such that $\dist_H(v, \phi(p_j))\le c$; 
			\item $w_1:E(H)\to \mathbb{N}$ is a map with size at most $c$, and $Q$ is a $w_1$-geodesic in $H$, and
			$\dist_{(H,w_1)}(r_i,r_j) = j-i$ for all $i,j\in J$ with $i<j$; and
			\item the subgraph of $H$ induced on the set of all $v\in V(H)$ with $\dist_H(v, \phi(P))>c$ belongs to $\mathcal{C}$, 
			where $\phi(P) = \{\phi(p_i):i\in I\}$. 
		\end{itemize}
		Then there is a function $w:E(H)\to \mathbb{N}$ with size at most $c_0$, such that $\phi$
		is a $(1,c_0)$-quasi-isometry from $G$ to $(H,w)$. 
	\end{thm}
	\Proof
	Let $r=2c(c+1)$, 
	and $c'=\max(\kappa(c),1)$.
	Let $c_2 = \max(2c+ c', 4(r+3)c^2)$.
	Define 
	$$c_3 = c_2+c(2(r+2)c+2)+  (r+2)cc' + (r+2)c^2,$$ 
	and 
	$$ c_0 = \max( (r+2)c^2, (r+2)cc', c_2, 4cc'+2c^2+2r+2crc_3).$$
	We will show that $c_0$ satisfies the theorem.
	
	Let $G,H, \phi, P$ and so on be as in the hypothesis of the theorem.
	Let $A$ be the set of all $v\in V(G)$ such that $\dist_G(v, P)\le r$. Let $B=V(G)\setminus A$.
	Let $X=\{\phi(v):v\in B\}$. 
	\\
	\\
	(1) {\em $\dist_H(X, \phi(P))\ge r/c-1$.}
	\\
	\\
	Let $b\in B$ and $i\in I$. 
	Then
	$$\dist_{H}(\phi(b), \phi(p_i))\ge (\dist_G(b,p_i)-c)/c\ge r/c-1.$$
	This proves (1).
	\\
	\\
	(2) {\em There is a partition $(Y,Z)$ of $V(H)\setminus X$, such that 
		\begin{itemize}
			\item for every $y\in Y$ there is a path of $H[X\cup Y]$ from $y$ to $X$, with length at most $(r+2)c$, and 
			$\dist_H(y,\phi(P))\ge  (r/c-1)/2>c$;
			\item for every $z\in Z$, there is a path of $H[Z]$ from $z$ to $\phi(P)$, with length at most $(r+2)c$, and 
			$\dist_H(z,X)> (r/c-1)/2$.
		\end{itemize}
	}
	\noindent
	Let $Y$ be the set of all $h\in V(H)\setminus X$ such that $\dist_H(h, X)\le \dist_H(h, \phi(P))$, and let 
	$Z=V(H)\setminus (X\cup Y)$.
	We claim that (2) is satisfied. Let $h\in V(H)\setminus X$. We claim first that either $\dist_H(h,X)\le c$, or 
	$\dist_H(h,\phi(P))\le cr+2c$. To see this, choose $v\in V(G)$ with $\dist_H(\phi(v), h)\le c$. If $v\in B$ then 
	$\phi(v)\in X$
	and the claim holds, so we assume that $v\in A$.
	Hence $\dist_G(v,P)\le r$, and so $\dist_H(\phi(v),\phi(P))\le cr+c$. Consequently
	$\dist_H(h, \phi(P))\le cr+2c$, and again the claim holds. Hence 
	$$\min(\dist_H(h, X), \dist_H(h, \phi(P)))\le (r+2)c,$$ 
	and so the first assertion of each bullet of (2) holds. For the second assertion, from (1), if $\dist_H(h, X)\le (r/c-1)/2$
	then $\dist_H(h, X)\le \dist(h,\phi(P))$ and therefore $h\in Y$; and similarly if $\dist_H(h, \phi(P))< (r/c-1)/2$ then 
	$h\in Z$. This proves (2).
	
	\bigskip
	
	Let $H'=H[X\cup Y]$. From (1) and (2), $\dist_H(y,\phi(P))> c$ for each $y\in X\cup Y$. Since
	the subgraph of $H$ induced on the set of all $v\in V(H)$ with $\dist_H(v, \phi(P))>c$ belongs to $\mathcal{C}$, by hypothesis, 
	and $\mathcal{C}$ is hereditary, it follows that $H'\in \mathcal{C}$.
	For each pair $b,b'\in B$, if $\dist_{H'}(\phi(b),\phi(b'))\le 2(r+2)c+1$, let $F_{b,b'}=F_{b',b}$ be a path between $b,b'$ of length
	$\dist_G(b,b')$, where all its internal vertices are new vertices. Let $F$ be the union of $G[B]$ and all the paths $F_{b,b'}$. 
	Define $\psi:V(F)\to V(H)$ as follows. For each $v\in B$, let $\psi(v) = \phi(v)$. For all $b,b'\in B$ and every internal vertex $v$ of 
	$F_{b,b'}$, let $\psi(v)$ be one of $\phi(b), \phi(b')$, chosen arbitrarily.
	\\
	\\
	(3) {\em If $u,v\in V(F)$, then $\dist_{H'}(\psi(u),\psi(v))\le (2(r+2)c+1)\dist_F(u,v)$.}
	\\
	\\
	It suffices to show that $\dist_{H'}(\psi(u),\psi(v))\le 2(r+2)c+1$ for every edge $uv$ of $F$ (and then sum over all
	edges of a geodesic of $F$ between $u,v$).
	Thus, let $uv\in E(F)$. If $uv$ is an edge of one of the paths $F_{b,b'}$, then 
	$$\dist_{H'}(\psi(u),\psi(v))\le \dist_{H'}(\phi(b), \phi(b'))\le 2(r+2)c+1,$$
	as required. 
	If $uv\in E(G[B])$, then 
	$\dist_H(\phi(u), \phi(v))\le 2c$
	since $\phi$ is a $(c,c)$-quasi-isometry from $G$ to $H$. Let 
	$S$ be a path
	of $H$ between $\phi(u),\phi(v)$ of length at most $2c$; so each of its vertices has distance at most $c$ from one of 
	$\phi(u),\phi(v)\in X$, and so $V(S)\subseteq X\cup Y$, since $c\le (r/c-1)/2$.
	Consequently, 
	$$\dist_{H'}(\psi(u),\psi(v))\le 2c\le 2(r+2)c+1.$$ 
	This proves (3).
	\\
	\\
	(4) {\em If $u,v\in V(F)$, then $\dist_F(u,v)\le 2c(2(r+2)c+1)\dist_{H'}(\psi(u),\psi(v))+4c(2(r+2)c+1)$.}
	\\
	\\
	Choose $u'\in B$ with $\psi(u) = \phi(u')$, and choose $v'$ similarly for $v$. 
	Let $T$ be a geodesic of $H'$ between $\phi(u'),\phi(v')$, and let its vertices be $t_0\LL t_n$ in order, where $t_0 = \phi(u')$
	and $t_n=\phi(v')$. For $0\le i\le n$, since $t_i\in X\cup Y$, there is a path $T_i$ of $H'$ from $t_i$ to $X$ with length at most $(r+2)c$; 
	let its end in $X$ be $x_i$, and choose $b_i\in B$ with $\phi(b_i) = x_i$. For $1\le i\le n$, there is a path from 
	$x_{i-1}$ to $x_{i}$ with vertex set a subset of 
	$V(T_{i-1})\cup V(T_{i})$, and its length is at most $2(r+2)c+1$; and consequently $F_{b_{i-1},b_{i}}$ exists, and so 
	$$\dist_F(b_{i-1},b_{i})=\dist_G(b_{i-1},b_{i})\le 2c\dist_H(x_{i-1},x_{i})\le 2c(2(r+2)c+1);$$
	so $\dist_F(b_{i-1},b_{i})\le 2c(2(r+2)c+1)$. But $\dist_F(b_0,b_n)$ is at most $\sum_{1\le i\le n}\dist_F(b_{i-1},b_{i})$
	and consequently  
	$$\dist_F(u', v')\le 2c(2(r+2)c+1) n =2c(2(r+2)c+1)\dist_{H'}(\psi(u),\psi(v)).$$
	But $\dist_F(u,u')\le 2c(2(r+2)c+1)$, and the same for $\dist_F(v,v')$;
	so 
	$$\dist_F(u,v)\le 2c(2(r+2)c+1)\dist_{H'}(\psi(u),\psi(v))+4c(2(r+2)c+1).$$
	This proves (4).
	
	\bigskip
	
	From the definition of $Y$, for each $y\in X\cup Y$ there exists $v\in V(F)$ such that $\dist_{H'}(\psi(v),y)\le (r+3)c$; and so,
	from (3) and (4), $\psi$
	is a $(2c(2(r+2)c+1), 4c(2(r+2)c+1))$-quasi-isometry from $F$ to $H'$. Since $\kappa$ is an additive bounder for $\mathcal{C}$, and $H'\in \mathcal{C}$,
	there is a function $w':E(H')\to \mathbb{N}$ with size at most $c'$, such that $\psi$
	is a $(1,c')$-quasi-isometry from $F$ to $(H',w')$,
	where $c'=\max(\kappa(c),1)$. Let $\Delta$ be the set of edges of $H$ between $X\cup Y$ and $Z$.
	Define $w:E(H)\to \mathbb{N}$ by:
	\begin{itemize}
		\item If $e\in E(H')$ then $w(e) = w'(e)$;
		\item If $e\in E(H[Z])$ then $w(e) = w_1(e)$;
		\item If $e\in \Delta$ then $w(e)=c_3$.
	\end{itemize}
	Thus $w$ has size at most $c_3$, and 
	we will show that $\phi$ is a $(1,c_0)$-quasi-isometry from $G$ to $(H,w)$.
	\\
	\\
	(5) {\em For all $i,j\in I$ with $i\le j$, $\dist_{(H,w)}(\phi(p_i),\phi(p_j))\le (j-i)+2c^2$.}
	\\
	\\
	Since for each $v\in V(Q)$ there exists $j\in J$ such that $\dist_H(v, \phi(p_j))\le c$, it follows that $V(Q)\subseteq Z$.
	From one of the hypotheses of the theorem, there exists $i'\in J$ with $|i'-i|\le c$ and $\dist_H(\phi(p_i), r_{i'})< c$;
	and there exists $j'\in J$ with $|j'-j|\le c$ and $\dist_H(\phi(p_j), r_{j'})< c$. Every geodesic of $H$ between 
	$\phi(p_i), r_{i'}$ has vertex set in $Z$, and so $\dist_{(H,w)}(\phi(p_i), r_{i'})\le (c-1)c$, since $w_1$ has size at most $c$;
	and similarly $\dist_{(H,w)}(\phi(p_j), r_{j'})\le (c-1)c$. Consequently $\dist_{(H,w)}(\phi(p_i), \phi(p_j))$ differs from 
	$\dist_{(H,w)}(r_{i'}, r_{j'})$ by at most $2(c-1)c$. 
	But $\dist_{(H,w)}(r_{i'},r_{j'})=|j'-i'|$, since $V(Q)\subseteq Z$;
	and so $\dist_{(H,w)}(\phi(p_i),\phi(p_j))\le |j'-i'|+2(c-1)c$. Since  $|i'-i|\le c$ and  $|j'-j|\le c$, it follows that 
	$|j'-i'|\le 2c+ (i'-i)$, and so 
	$$\dist_{(H,w)}(\phi(p_i),\phi(p_j))\le (i'-i) + 2(c-1)c+2c.$$
	This proves (5).
	\\
	\\
	(6) {\em Let $u,v\in V(G)$. Then 
		$$\dist_{(H,w)}(\phi(u), \phi(v))\le \dist_G(u,v)+ 4cc' + 2c^2 + 2r+ 2crc_3.$$}
	\noindent
	Observe first that if $T$ is a geodesic of $G$, with $V(T)\subseteq B$ and with ends $b_1,b_2$ say, then 
	$$\dist_{(H,w)}(\phi(b_1), \phi(b_2))\le \dist_{(H',w')}(\psi(b_1), \psi(b_2)) \le \dist_F(b_1,b_2)+c'=\dist_G(b_1,b_2)+c',$$ 
	from the choice of $w'$. Now let $T$ be a geodesic in $G$ between $u,v$; and we may therefore assume that
	$V(T)\not\subseteq B$. Let $a_1,a_2$ be the first and last vertices of $T$ that belong to $A$. If $a_1\ne u$, let $b_1\in V(T)$ be 
	adjacent in $T$ to $a_1$, and not between $a_1,a_2$; thus $b_1\in B$ from the definition of $a_1$. Let $T_1=T[u,b_1]$. 
	If $a_1=u$ then $b_1,T_1$
	are undefined. Define $b_2,T_2$ similarly if $a_2\ne v$. 
	
	If $b_1,T_1$ exist, then $T_1$ is a geodesic of $G$ with vertex set in $B$, and so 
	$$\dist_{(H,w)}(\phi(u), \phi(b_1))\le \dist_G(u,b_1)+c',$$ 
	as we showed above.
	Since $a_1b_1\in E(G)$ and $\phi$ is a $(c-1,c)$-quasi-isometry from $G$ to $H$, it follows that $\dist_H(\phi(a_1),\phi(b_1))\le 2c-1$.
	Consequently $\dist_{H'}(\phi(a_1),\phi(b_1))\le 2c-1$, as the corresponding path in $H$ is
	contained in $H'$;
	and since $w'$ has size at most $c'$,
	it follows that $\dist_{(H,w)}(\phi(a_1),\phi(b_1))\le (2c-1) c'$.
	Thus, if $b_1,T_1$ exist, then 
	$$\dist_{(H,w)}(\phi(u), \phi(a_1))\le 
	\dist_G(u,b_1)+c'+ (2c-1) c' \le \dist_G(u,a_1)+2c c'.$$
	This last is also trivially true if $b_1,T_1$ do not exist, since then $u=a_1$. A similar inequality holds for $v,b_2$.
	
	Since $a_1\in A$, there exists $i_1\in I$ such that 
	$\dist_G(a_1,p_{i_1})\le r$. Choose $i_2$ similarly for $a_2$. Thus $\dist_G(p_{i_1},p_{i_2})\le \dist_G(a_1,a_2)+2r$, and so 
	$\dist_G(a_1,a_2)\ge |i_2-i_1|-2r$. Now since $\dist_G(a_1,p_{i_1})\le r$,
	and $\phi$ is a $(c-1,c)$-quasi-isometry from $G$ to $H$, it follows that $\dist_H(\phi(a_1), \phi(p_{i_1}))\le (c-1)r+c\le cr$, and 
	so $\dist_{(H,w)}(\phi(a_1), \phi(p_{i_1}))\le cr c_3$. The same holds for $a_2, p_{i_2}$; and so 
	$$\dist_{(H,w)}(\phi(a_1), \phi(a_2))\le \dist_{(H,w)}(\phi(p_{i_1}), \phi(p_{i_2}))+ 2cr c_3.$$
	Since 
	$$\dist_{(H,w)}(\phi(p_{i_1}), \phi(p_{i_2}))\le |i_2-i_1|+2c^2$$ 
	by (5), we deduce that
	$$\dist_{(H,w)}(\phi(a_1), \phi(a_2))\le  |i_2-i_1|+2c^2+ 2cr c_3.$$
	But $\dist_G(a_1,a_2)\ge |i_2-i_1|-2r$, and so 
	$$\dist_{(H,w)}(\phi(a_1), \phi(a_2))\le \dist_G(a_1,a_2) +2r + 2c^2+ 2crc_3.$$
	We deduce that 
	\begin{align*}
		\dist_{(H,w)}&(\phi(u),\phi(v))\\
		&\le \dist_{(H,w)}(\phi(u), \phi(a_1)) + \dist_{(H,w)}(\phi(a_1), \phi(a_2)) + \dist_{(H,w)}(\phi(v), \phi(a_2))\\
		&\le \dist_G(u,a_1)+2cc' + \dist_G(a_1,a_2) +2r+2c^2+ 2cr c_3 + \dist_G(v,a_2)+2cc'\\
		&=\dist_G(u,v) + 4cc' + 2c^2 +2r + 2crc_3 .
	\end{align*}
	This proves (6).
	\\
	\\
	(7) {\em Let $T$ be a path of $H[Z]$ between $\phi(a_1), \phi(a_2)$, where $a_1,a_2\in V(G)$.
		Then 
		$$\dist_G(a_1,a_2)\le w(T)+ 4(r+3)c^2.$$}
	\noindent
	For $t = 1,2$, since $\phi(a_t)\in Z$, there exists $i_t\in I$ such that there is a path of $H[Z]$ between 
	$\phi(a_t), \phi(p_{i_t})$ of length at most $(r+2)c$, and 
	there exists $j_t\in J$ such that 
	$|j_t-i_t|\le c$ and $\dist_H(\phi(p_{i_t}), r_{j_t})< c$. 
	Since $\phi$ is a $(c-1,c)$-quasi-isometry, it follows that 
	$\dist_G(a_t, p_{i_t})\le (c-1)(r+2)c+c\le (r+2)c^2$ for $t=1,2$, and so 
	$$\dist_G(a_1,a_2)\le \dist_G(p_{i_1},p_{i_2})+2c^2(r+2).$$
	But 
	$$\dist_G(p_{i_1},p_{i_2}) = |i_2-i_1|\le |j_2-j_1|+2c=\dist_{(H,w_1)}(r_{j_1},r_{j_2})+2c,$$
	and so 
	$$\dist_G(a_1,a_2)\le \dist_{(H,w_1)}(r_{j_1},r_{j_2})+2(r+2)c^2+2c.$$
	Since for $t = 1,2$ there is a path of $H$ between $\phi(a_t),r_{j_t}$ of length at most $(r+2)c+c$, and hence 
	$\dist_{(H,w_1)}(\phi(a_t),r_{j_t})\le (r+3)c^2$, and $w_1(T)=w(T)$, it follows that 
	$$\dist_{(H,w_1)}(r_{j_1},r_{j_2})\le \dist_{(H,w_1)}(\phi(a_1),\phi(a_2))+2(r+3)c^2\le w(T)+2(r+3)c^2.$$
	Thus 
	$$\dist_G(a_1,a_2)\le w(T)+2(r+3)c^2+2(r+2)c^2+2c\le w(T)+4(r+3)c^2.$$
	This proves (7).
	\\
	\\
	(8) {\em Let $u,v\in V(G)$, and let $T$ be a path of $H$ between $\phi(u), \phi(v)$. 
		Then $\dist_G(u,v)\le w(T) + c_2$.}
	\\
	\\
	We proceed by induction on $|\Delta\cap E(T)|$. Suppose first that $\Delta\cap E(T)=\emptyset$, and so $T$ is a path of one of $H'$, $H[Z]$. 
	If $T$ is a path of $H[Z]$, the result holds by (7), so 
	we assume that $T$ is a 
	path of $H'$. Thus there exist $b_1,b_2\in B$ with $\phi(b_1) = \phi(u)$ and $\phi(b_2) = \phi(v)$. Since $\phi$ is a $(c,c)$-quasi-isometry
	from $G$ to $H$, it follows that $\dist_G(u,b_1), \dist_G(v, b_2)\le c$. Moreover, $\dist_{H'}(\phi(u), \phi(v))\le w(T)$, and so 
	$\dist_G(b_1,b_2)\le \dist_{F}(b_1,b_2)\le w(T)+c'$, since $\phi$ is a $(1,c')$-quasi-isometry
	from $H'$ to $H[X\cup Y]$. It follows that in this case, $\dist_G(u,v)\le w(T)+ 2c+ c'$, and so the result holds.

	Thus we may assume that there exists $yz\in \Delta\cap E(T)$, where $y\in X\cup Y$ and $z\in Z$. 
	By exchanging $u,v$ if necessary we may 
	assume that $\phi(u), y,z,\phi(v)$ are in order in $T$. Since $y\in X\cup Y$, there exists $b\in B$ such that 
	$\dist_{H'}(\phi(b), y)\le (r+2)c$, and hence $\dist_{(H,w)}(\phi(b), y)\le (r+2)cc'$, since $w'$ has size at most $c'$;
	and since $z\in Z$, there exists $i\in I$ such that $\dist_{H[Z]}(z, \phi(p_i))\le (r+2)c$, and hence 
	$\dist_{(H,w)}(z, \phi(p_i))\le (r+2)c^2$. By combining these paths with the subpaths $T[\phi(u), y]$ and $T[z, \phi(v)]$ respectively,
	we deduce (since $w(yz) = c_3$) that there are paths $R_1,R_2$ of $H$, where $R_1$ is between $\phi(u), \phi(b)$, and $R_2$ is between $\phi(p_i), \phi(v)$, such that
	$$w(R_1)+w(R_2)\le w(T)-c_3+(r+2)cc'+(r+2)c^2,$$
	and $R_1,R_2$ both have fewer than $|\Delta\cap E(T)|$ edges in $\Delta$.
	From the inductive hypothesis,
	$\dist_G(u,b)\le w(R_1) +c_2$, and $\dist_G(p_i,v)\le w(R_2) + c_2$.
	But 
	$$\dist_G(u,v)\le \dist_G(u,b) + \dist_G(b,p_i) +\dist_G(p_i,v),$$
	and 
	$$\dist_G(b,p_i)\le c\dist_H(\phi(b), \phi(p_i))+c\le c(2(r+2)c+1)+c;$$
	so 
	\begin{align*}
		\dist_G(u,v) &\le \dist_G(u,b) + \dist_G(p_i,v) + c(2(r+2)c+2)\\
		&\le w(R_1) + c_2 + w(R_2) +c_2 + c(2(r+2)c+2)\\
		&\le  w(T)+ 2c_2 + c(2(r+2)c+2) -c_3+  (r+2)cc' + (r+2)c^2\\
		&\le w(T)+c_2.
	\end{align*}
	This proves (8).
	\\
	\\
	(9) {\em For each $v\in V(H)$, there exists $u\in V(G)$ such that $\dist_{(H,w)}(\phi(u), v)\le (r+2)c\max(c,c')\le c_0$.}
	\\
	\\
	If $v\in Z$, then by (2), there is a path of $H[Z]$ from $v$ to $\phi(P)$, of length at most $(r+2)c$, and hence
	$\dist_{(H,w)}(\phi(u), v)\le (r+2)c^2$. If $v\in X\cup Y$, by (2) there is a path of $H'$ from $v$ to $X$, of length at most $(r+2)c$, 
	and hence $\dist_{(H,w)}(v,X)\le (r+2)cc'$. This proves (9). 
	
	\bigskip
	
	By (6), (8) and (9),  $\phi$ is a $(1,c_0)$-quasi-isometry from $G$ to $(H,w)$, and its size is $c_3$. 
	This proves \ref{usegeo}.~\bbox
	
	
	\section{Combining the two steps}
	
	In this section we complete the proof of \ref{weightedline} when $H$ is finite, and reduce the problem to finding the geodesic $P$
	when $H$ is infinite. 
	If $\phi$ is a quasi-isometry from a graph $G$ to $H$, and $X\subseteq V(H)$, we denote by $\phi^{-1}(X)$ the set
	of all $v\in V(G)$ with $\phi(v)\in X$.
	We would like to prove \ref{weightedline} by induction on the line-width of $H$, and the next result is half of the inductive step.
	\begin{thm}\label{halfinduct}
		Let $k\ge 1$ be an integer, and suppose that for all $C\ge 1$ there exists $C'$ with the following property.
		\begin{itemize}
			\item If $\phi$ is a $(C-1,C)$-quasi-isometry from a graph $G$ to a graph $H$ with
			line-width at most $k-1$, then there is a function $w:E(H)\to \mathbb{N}$ with size at most $C'$, such that $\phi$
			is a $(1,C')$-quasi-isometry from $G$ to $(H,w)$. 
		\end{itemize}
		Then for all $C,D\ge 1$ there exists $C'\ge 0$ with the following property. 
		\begin{itemize}
			\item Suppose that 
			$\phi$ is a $(C-1,C)$-quasi-isometry from a graph $G$ to a graph $H$ with
			a line-decomposition $(B_t:t\in T)$ of width at most $k$; and $P$ is a geodesic of $G$ 
			such that 
			$\dist_{G}(P,\phi^{-1}(B_t))\le D$ for each $t\in V(T)$.
			Then there is a function $w:E(H)\to \mathbb{N}$ with size at most $C'$, such that $\phi$
			is a $(1,C')$-quasi-isometry from $G$ to $(H,w)$. 
		\end{itemize}
	\end{thm}
	\Proof Let $C,D\ge 1$; and we may assume that $C\ge 2$. From the hypothesis, there is an additive bounder $\kappa$ for the class 
	$\mathcal{C}$ of all graphs with line-width less than $k$.
	Let $c_0$ be as in \ref{usegeo}, taking 
	$c=\max(32C^4, CD+C)$.
	Define $C'=c_0$; we will show that $C'$ satisfies the theorem. 
	
	Let $\phi$ be a $(C-1,C)$-quasi-isometry from a graph $G$ to a graph $H$ with
	a line-decomposition $(B_t:t\in T)$ of width at most $k$, and let $P$ be 
	a geodesic of $G$ such that 
	$\dist_{G}(P,\phi^{-1}(B_t))\le D$ for each $t\in T$.
	Let the vertices of $P$ be $p_i\;(i\in I)$ in order, where $I$ is an integer interval.
	By \ref{fixgeo}, 
	there is a function $w_1:E(H)\to \mathbb{N}$, with size at most $32C^4$, and a path $Q$ of $H$, and $J\subseteq I$, cofinal with $I$, and distinct vertices
	$r_j\;(j\in J)$, in order in $Q$, with the following properties:
	\begin{itemize}
		\item $Q$ is a $w_1$-geodesic in $H$;
		\item $\dist_{(H,w_1)}(r_i,r_j) = j-i$ for all $i,j\in J$ with $i<j$;
		\item for all $i\in I$ there exists $j\in J$ with $|j-i|\le C^2$ and $\dist_H(\phi(p_i), r_j)< C^3$;
		\item $\dist_H(\phi(p_j), r_j)\le 2C$ for each $j\in J$; and
		\item for each $v\in V(Q)$ there exists $j\in J$ such that $\dist_H(v, \phi(p_j))\le C$.
	\end{itemize}
	We may assume that $Q$ is the union of the subpaths $Q[r_i,r_j]$ for $i,j\in J$, by replacing $Q$ by this union if necessary.
	\\
	\\
	(1) {\em For each $t\in T$, there exists $h\in B_t$ such that $\dist_H(h,\phi(P))\le  CD+C$.}
	\\
	\\
	By hypothesis, $\dist_{G}(P,\phi^{-1}(B_t))\le D$; choose $v\in \phi^{-1}(B_t)$ with $\dist_{G}(P,v)\le D$, and let $h=\phi(v)$. 
	Hence 
	$\dist_H(\phi(v),\phi(P))\le (C-1)D+C\le CD+C$. This proves (1).
	
	\bigskip
	
	From (1), the subgraph of $H$ induced on the set of all $h\in V(H)$ with $\dist_H(h, \phi(P))>CD+c$ has line-width at most $k-1$,
	where $\phi(P) = \{\phi(p_i):i\in I\}$. By \ref{usegeo}, taking $c=\max(32C^4, CD+C)$, we deduce that 
	there is a function $w:E(H)\to \mathbb{N}$ with size at most $c_0$, such that $\phi$
	is a $(1,c_0)$-quasi-isometry from $G$ to $(H,w)$. This proves \ref{halfinduct}.~\bbox
	
	To complete the proof of \ref{weightedline} by induction on the line-width of $H$, it therefore suffices to obtain a geodesic $P$ with
	the properties of \ref{halfinduct}. This is simple in the finite case, and lengthy in the infinite case, so let us do the finite case separately, for readers whose only interest is the finite case. We need: 
	\begin{thm}\label{shortjump}
		Let $(B_t:t\in T)$ be a line-decomposition of a graph $H$, and let $C\ge 1$.
		Let $\phi$ be a surjective $(C-1,C)$-quasi-isometry from a graph $G$ to $H$. Let $t,t',t''\in T$ with $t'\le t\le t''$, and
		let $K$ be a connected subgraph of $G$
		with $V(K)\cap \phi^{-1}(B_{t'}), V(K)\cap \phi^{-1}(B_{t''})$ both nonempty.
		Then there is a vertex $x\in V(K)$ such that $\dist_H(\phi(x),B_t)\le C-1$.
	\end{thm}
	\Proof
	Since $V(K)\cap \phi^{-1}(B_{t'}), V(K)\cap \phi^{-1}(B_{t''})$ are both nonempty and $K$ is connected, there is a path $P$ of $K$
	with ends $x', x''$ say, where $\phi(x')\in B_{t'}$ and $\phi(x'')\in B_{t''}$.
	If $\phi(x)\in B_{t}$ for some $x\in V(P)$, then $x$ satisfies the theorem, so we suppose not.
	Since every path in $H$ between $B_{t'}$ and $B_{t''}$ has a vertex in $B_{t}$, it follows that
	$\phi(x'), \phi(x'')$ belong to different components of $H\setminus B_{t}$. Consequently there is an edge $ab$ of $P$
	such that $\phi(a),\phi(b)$ belong to different components of $H\setminus B_{t}$. Since $\dist_H(\phi(a), \phi(b))\le 2C-1$,
	there exists $x\in \{a,b\}$ such that $\dist_H(\phi(x), B_t)\le C-1$.
	This proves \ref{shortjump}.~\bbox
	
	Let $H$ be a graph that admits a line-decomposition $(B_t:t\in T)$. We can remove from $T$
	all $t$ with $B_t=\emptyset$, so we may assume that $B_t\ne \emptyset$ for each $t\in T$, that is, $(B_t:t\in T)$ is {\em nowhere-null}.
	
	\begin{thm}\label{finiteend}
		Suppose that
		$\phi$ is a $(C-1,C)$-quasi-isometry from a connected graph $G$ to a finite graph $H$ with
		a nowhere-null line-decomposition $(B_t:t\in T)$ of width at most $k$. Then there is a geodesic $P$ of $G$
		such that
		$\dist_{G}(P,\phi^{-1}(B_t))\le C^2$ for each $t\in V(T)$.
	\end{thm}
	\Proof
	Since $H$ is finite, we may assume that $T$ is finite, and so we may assume that $T=\{1\LL n\}$. Choose $a\in B_1$ and $b\in B_n$, 
	and choose $u,v\in V(G)$ such that $\dist_H(\phi(u),a)\le C$, and $\dist_H(\phi(v),b)\le C$. Let $P$ be a geodesic in $G$ between $u,v$
	(this exists since $G$ is connected). Choose $h\in \{1\LL n\}$ minimum such that $\phi(u)\in B_h$, and $j\in \{1\LL n\}$ maximum such that
	$\phi(v)\in B_j$. Now let $t\in \{1\LL n\}$; we need to show that $\dist_{G}(P,\phi^{-1}(B_t))\le C^2$. If $h\le t\le j$, then 
	this is true by \ref{shortjump}, since $C^2\ge C$; so from the symmetry we may assume that $t<h$. There is a path of $H$ between $B_1$ and
	$\phi(u)$ of length at most $C$, and so some vertex of this path belongs to $B_t$; and consequently $\dist_H(\phi(u),B_t)\le C$, and so
	$\dist_G(u,\phi^{-1}(B_t))\le C^2$.  This proves \ref{finiteend}.~\bbox

	Thus, to complete the proof of \ref{weightedline} when $H$ is finite, we may assume that $H$ is connected, and hence $G$ is connected
	(because there is a quasi-isometry between them). Choose $P$ as in \ref{finiteend}; then it satisfies the hypothesis of \ref{halfinduct}, and so \ref{halfinduct}  completes the inductive proof.

	\section{Some simplifications}
	
	The remainder of the paper concerns obtaining the geodesic $P$ when $H$ is infinite. As before, we can assume that $G,H$ are connected.
	Before the main argument, in the next section, we take off some bite-sized pieces.
	Let us see first that, to prove \ref{weightedline} in general, it suffices to prove the result when $\phi$ is surjective.
	
	\begin{thm}\label{reducetosur}
		Let $L,C\ge 0$; and suppose that there exist $C', W$ such that if $\phi$ is a surjective $(\lfloor |L/(2C+1)\rfloor,C)$-quasi-isometry from a graph $G$ to a graph $H$ with
		line-width at most $k$, then there is a function $w:E(H)\to \mathbb{N}$ such that $\phi$
		is a $(1,C')$-quasi-isometry from $G$ to  $(H,w)$. It follows that, if $\phi$ is an $(L,C)$-quasi-isometry from a graph $G$ to a graph $H$ with
		line-width at most $k$, then there is a function $w:E(H)\to \mathbb{N}$ such that $\phi$
		is a $(1,C')$-quasi-isometry from $G$ to  the weighted graph $(H,w)$.
	\end{thm}
	\Proof 
	Let $L'=\lfloor |L/(2C+1)\rfloor$, and suppose that  $\phi$ is an $(L,C)$-quasi-isometry from a graph $G$ to a graph $H$ with
	line-width at most $k$. Let $Z=\{\phi(v):v\in V(G)\}$.
	For each $x\in V(H)\setminus Z$, since
	$\phi$
	is an $(L,C)$-quasi-isometry, there exists $z\in Z$ such that $\dist_H(z,x)\le C$. We deduce (by adding a new vertex $r$ adjacent to each vertex in $Z$,
	choosing a breadth-first tree rooted at $r$, and then deleting $r$) that for each $z\in Z$, there exists a subset 
	$\eta(z)\subseteq V(H)$ with the following properties:
	\begin{itemize}
		\item $\eta(z) \cap Z=\{z\}$, and the sets $\eta(z)\;(z\in Z)$ are pairwise disjoint and have union $V(H)$;
		\item for each $z\in Z$ and each $v\in \eta(z)$, there is a path $P$ of $H[\eta(z)]$ between $v,z$ with length at most $C$.
	\end{itemize}
	Let $H_1$ be obtained from $H$ by contracting each of the connected subgraphs $H[\eta(z)]$ to a single vertex $p(z)$.
	It follows that $H_1$ has line-width at most $k$. For each $v\in V(G)$, let $\phi_1(v) = p(\phi(v))$.
	It follows that if $u,v\in V(G)$, then
	$$\dist_{H_1}(\phi_1(u),\phi_1(v))\le \dist_{H}(\phi(u), \phi(v))\le L'\dist_G(u,v)+C.$$
	Also, let $P$ be a geodesic of $H_1$ between $\phi_1(u),\phi_1(v)$, with vertices $p(z_0)\CC p(z_n)$ in order, say, where
	$z_0\LL z_k\in Z$. So $z_0=\phi(u)$ and $z_n=\phi(v)$.
	For $0\le i<k$, let $e_i$ be the edge of $H_1$ between $p(z_i), p(z_{i+1})$. Then $e_i$ is an edge of $H$ between
	$\eta(z_i),\eta(z_{i+1})$, and so there is a path $P_i$ of $H$ between $z_i, z_{i+1}$ with length at most $1+ 2C$.
	By concatenating these paths, we find that
	$$\dist_{H}(z_0, z_n)\le (2C+1)|E(P)|=(2C+1)\dist_{H_1}(\phi_1(u),\phi_1(v)).$$
	Consequently
	$$\dist_G(u,v)\le L'\dist_{H}(\phi(u),\phi(v))+C \le L\dist_{H_1}(\phi_1(u),\phi_1(v))+C.$$
	It follows that $\phi_1$ is a surjective $(L, C)$-quasi-isometry from $G$ to $H_1$, and $H_1$ is connected.
	
	Thus, there is a function $w_1: E(H_1)\to \mathbb{N}$ with size at most $W$ such that $\phi_1$ is a $(1,C')$-quasi-isometry
	from $G$ to $(H_1,w_1)$. Let $w(e) = w_1(e)$ for each $e\in E(H_1)$, and let $w_1(e) = 0$ for each edge $e$ of $H$ that is not an
	edge of $H_1$ (and so has both ends in $\eta(z)$ for some $z\in Z$). It follows that for all $z,z'\in Z$,
	$$\dist_{(H,w')}(z,z') = \dist_{H_1}(p(z), p(z')), $$
	and consequently, for all $u,v\in V(G)$,
	$$\dist_{(H,w)}(\phi(u), \phi(v)) = \dist_{H_1}(\phi_1(u), \phi_1(v)).$$
	Since $\phi_1$ is a $(1,C')$-quasi-isometry
	from $G$ to $(H_1,w_1)$, and $\dist_{(H,w)}(x,Z)=0$ for each $x\in V(H)\setminus Z$, it follows that $\phi$ is a $(1,C')$-quasi-isometry
	from $G$ to $(H,w)$.
	This proves \ref{reducetosur}.~\bbox

	We will prove the following in the next section:
	\begin{thm}\label{addgeo0}
		Let $(B_t:t\in T)$ be a nowhere-null line-decomposition of a connected graph $H$, with finite width, and let $C\ge 2$.
		Let $\phi$ be a surjective $(C-1,C)$-quasi-isometry from a graph $G$ to $H$.
		Then there is a graph $G'$ and a geodesic $P$ of $G'$ with the following properties:
		\begin{itemize}
			\item $G$ is an induced subgraph of $G'$, and $\dist_G(u,v) = \dist_{G'}(u,v)$ for all $u,v\in V(G)$;
			\item the identity map from $V(G)$ into $V(G')$ is a $(1,2C^2+1)$-quasi-isometry, and there is a
			$(1,4C^2+2)$-quasi-isometry from $G'$ to $G$ that maps each vertex of $G$ to itself; and
			\item for each $t\in T$, $\dist_{G'}(P,\phi^{-1}(B_t))\le 6C^2$.
		\end{itemize}
	\end{thm}
	
	Next, we show that if \ref{addgeo0} holds then we can complete the proof of \ref{weightedline}.
	We need:
	\begin{thm}\label{compose}
		Let $\phi$ be an $(L_1,C_1)$-quasi-isometry from $G$ to $H$, and let $\psi$ be an $(L_2, C_2)$-quasi-isometry from $F$ to $G$.
		For each $v\in V(F)$, define $\theta(v) = \phi(\psi(v))$. Then $\theta$ is an $(L_1L_2,C)$-quasi-isometry from $F$ to $H$,
		provided that $C\ge \max (L_1C_2+2C_1, L_2C_1+C_2)$.
	\end{thm}
	The proof is routine calculation and we omit it.
	
	\bigskip
	
	\noindent {\bf Proof of \ref{weightedline}, assuming \ref{addgeo0}:\ \ }
	We proceed by induction on $k$, and we already saw in \ref{reducetosur} that it suffices to prove \ref{weightedline} when $\phi$ is surjective.
	Let $L,C\ge 0$. Let $C_0=L(2C^2+2)+2C$. Choose $C'$ as in \ref{halfinduct}, with $C,D$ replaced by $C_0, 6C^2$ respectively, and 
	let $W=C'$. We claim that $C', W$ satisfy \ref{weightedline}, when $\phi$ is surjective. 
	
	Let $H$ be a graph with line-width at most $k$
	and let $\phi$ be a surjective $(L,C)$-quasi-isometry from a connected graph $G$ to $H$.
	By \ref{addgeo0},
	there is a graph $G'$ and a geodesic $P$ of $G'$ with the following properties:
	\begin{itemize}
		\item $G$ is an induced subgraph of $G'$, and $\dist_G(u,v) = \dist_{G'}(u,v)$ for all $u,v\in V(G)$;
		\item the identity map from $V(G)$ into $V(G')$ is a $(1,2C^2+1)$-quasi-isometry, and there is a
		$(1,4C^2+2)$-quasi-isometry $\psi$ from $G'$ to $G$ that maps each vertex of $G$ to itself; and
		\item for each $t\in T$, $\dist_{G'}(P,\phi^{-1}(B_t))\le 6C^2$.
	\end{itemize}
	
	For each $v\in V(G')$, define $\theta(v) = \phi(\psi(v))$. Since this is the composition of an $(L,C)$-quasi-isometry and a
	$(1,4C^2+2)$-quasi-isometry, it follows from \ref{compose} that $\theta$ is an $(L, L(2C^2+2)+2C)$-quasi-isometry from $G'$ to $H$,
	and hence a $(C_0-1, C_0)$-quasi-isometry.
	Moreover, 
	$$\dist_{G'}(P,\phi^{-1}(B_t))\le 6C^2$$ 
	for each $t\in T$, and so we can apply \ref{halfinduct}. 
	We deduce that there is a function $w:E(H)\to \mathbb{N}$ with size at most $C'$, such that $\theta$
	is a $(1,C')$-quasi-isometry from $G'$ to $(H,w)$. 
	Since $\phi$ is surjective, and
	$\dist_G(u,v) = \dist_{G'}(u,v)$ for all $u,v\in V(G)$, it follows that the restriction of $\theta$ to $V(G)$ is also a
	$(1,C')$-quasi-isometry from $G$ to $(H,w)$. But this restriction is just $\phi$, since $\psi$ maps each vertex of $G$ to itself.
	This proves \ref{weightedline}.~\bbox

	\section{Finding a spanning geodesic}
	
	We begin with an example.  Construct graphs $G, H$ as follows.
	For each $j\ge 1$, let $Q_j$ be a path of length $2j$, all vertex-disjoint,
	and for each $j$ let $Q_{j}$ have vertices
	$$q_{j}^{-j}\DD q_{j}^{-j+1}\CC q_{j}^{0}\CC q_{j}^{j-1}\DD q_{j}^{j},$$
	in order.
	For each integer $i$ (including negative integers) let $v_i$ be a new vertex,
	adjacent to $q_j^i$ for each $j$ with $1\le j\le |i|$, forming $G$.
	Every geodesic of $G$ is finite, because
	every geodesic contains at most two of the vertices $v_{i}\;(i\in \mathbb{Z})$.
	
	Let $H$ be the two-way infinite path with vertex set $\{v_{i}:i\in \mathbb{Z}\}$,  where $v_i,v_{i+1}$ are adjacent for each $i$.
	For each $v\in V(G)$, let $\phi(v)=v$
	if $v\in V(H)$, and let $\phi(v) = v_i$ if $v=q_j^i$. Then $\phi$ is a surjective $(1,1)$-quasi-isometry from $G$ to $H$. For each integer $i$ let $B_i=\{v_i, v_{i+1}\}$; so $(B_i:i\in \mathbb{Z})$ is a path-decomposition of $H$.
	It follows easily that there is no geodesic $P$ in $G$ (with respect to this path-decomposition of $H$) as we need for \ref{halfinduct}.
	
	Thus $P$ might not always exist, and this section concerns faking
	up a substitute. 
	An {\em interval} of a linearly ordered set $T$ is a set $J\subseteq T$
	such that if $t_1,t_2,t_3\in T$ with $t_1<t_2<t_3$ and $t_1,t_3\in J$, then $t_2\in J$.
	If $(B_t:t\in T)$ is a line-decomposition of $H$, then for each $v\in V(H)$, the set $\{t\in T:v\in B_t\}$ is an interval of $T$,
	and we denote it by $\tau(v)$.
	
	If $T$ is a totally ordered set, we say that $J\subseteq T$ is an {\em initial interval} of $T$ if there do not exist
	$i,j\in T$ with $i\le j$ and $j\in J$ and $i\notin J$; and $J\subseteq T$ is a {\em final interval} of $L$ if there do not exist
	$i,j\in T$ with $i\le j$ and $i\in J$ and $j\notin J$.
	If $L\subseteq T$ is a final interval of $T$, we say that $v\in V(H)$ is an {\em southern border vertex} (for $L$)
	if $\tau(v)\cap L\ne \emptyset$ and $\tau(v)\not\subseteq L$.
	
	We begin with:
	\begin{thm}\label{connected}
		Let $H$ be a connected graph, and let $(B_t:t\in T)$ be a nowhere-null line-decomposition of $H$ with finite width.
		Let $J$ be an initial interval of $T$ and let $L=T\setminus J$.
		If $J, L\ne \emptyset$, then there is a southern border vertex for $L$, and there exists $t\in L$ containing all
		such vertices.
	\end{thm}
	\Proof Let $Y$ be the set of southern border vertices for $L$.
	Let $P=\bigcup_{t\in J}B_t$ and $Q = \bigcup_{t\in L}B_t$. Thus $P\cup Q=V(H)$, and we suppose (for a contradiction) that
	$P\cap Q=\emptyset$.
	Since $(B_t:t\in T)$ is non-null, and $J,L\ne \emptyset$, it follows that $P,Q\ne \emptyset$. Since $H$ is connected,
	there is an edge of $H$ joining $P,Q$, and so there exists $t\in T$ such that $B_t\cap P,B_t\cap Q\ne \emptyset$.
	But then $t\notin J\cup L$,
	a contradiction. This proves that $P\cap Q\ne \emptyset$. Since $P\cap Q\subseteq Y$,
	we deduce that $Y\ne \emptyset$.
	
	We claim that the intervals $\tau(y)\;(y\in Y)$ pairwise intersect. To see this, let $y_1,y_2\in Y$, and choose
	$\ell_1\in L\cap \tau(y_1)$ and $\ell_2\in L\cap \tau(y_2)$. We may assume that $\ell_1\le \ell_2$; but then $\ell_1\in \tau(y_2)$
	since $y_2$ is a southern border vertex of $L$. This proves that $\tau(y_1),\tau(y_2)$ intersect, and so
	the intervals $\tau(y)\;(y\in Y)$ pairwise intersect.
	
	From the finite Helly property of intervals, it follows that for every finite subset $Y'$ of $Y$, there exists $t\in L$
	that belongs to $\tau(y)$
	for each $y\in Y'$. Since $(B_t:t\in T)$ has finite width $k$ say, it follows that $|Y'|\le |B_t|\le k+1$, and so $Y$ is finite, and
	therefore there exists $t\in L$ with $Y\subseteq B_t$.
	This proves \ref{connected}.~\bbox

	If $(B_t:t\in T)$ is a non-null line-decomposition of a graph $H$,
	we say that $Y\subseteq T$ is {\em upfinal} in $T$ if
	\begin{itemize}
		\item for each $t\in T$, there exists $y\in Y$ such that $y\ge t$; and
		\item for each $t\in T$, there are only finitely many $y\in Y$ with $y<t$.
	\end{itemize}
	We define {\em downfinal} similarly.
	We say that $X\subseteq V(H)$ is {\em up-pervasive} if
	\begin{itemize}
		\item for each $t\in T$, there exists $t'\in T$ with $t'\ge t$ such that $X\cap B_{t'}\ne \emptyset$; and
		\item for each $t\in T$, only finitely many vertices of $X$ belong to $\bigcup_{t'\le t}B_{t'}$.
	\end{itemize}
	{\em Down-pervasive} is defined similarly, reversing the order of $T$.
	
	There might be a vertex $x$ such that $\{x\}$ is up-pervasive; that is, $\tau(x)$ is a final interval of $T$. Similarly there might be
	$x'$ such that $\{x'\}$ is down-pervasive. If both exist then the problem of this section is easy to resolve, and if even one exists,
	it helps a good deal. The main case is when neither exists. Let us say that $(B_t:t\in T)$ is {\em upper-open} if no
	singleton set is up-pervasive. Every finite up-pervasive set includes a singleon up-pervasive set, so if  $(B_t:t\in T)$ is upper-open
	then no finite set is up-pervasive. {\em Lower-open} is defined similarly. 
	
	There need not be an finite up-pervasive set, but for the line-decompositions of concern to us, there is always one that is countable, as the next result shows.
	
	
	\begin{thm}\label{getupfinal}
		Let $(B_t:t\in T)$ be a non-null upper-open line-decomposition of a connected graph $H$, with finite width, and let $t_0\in T$.
		Then there is an infinite sequence $t_1,t_2\LL$ of elements of $T$, 
		with the following properties:
		\begin{itemize}
			\item $t_i<t_j$ for all $i,j$ with $0\le i<j$;
			\item for each $t\in T$ with $t\ge t_0$, there is at least one and at most two values of $i\ge 0$ such that $B_t\cap B_{t_i}\ne \emptyset$; and
			\item for each $t\in T$, there are only finitely many $i\ge 0$ such that $t_i\le t$.
		\end{itemize}
		Consequently there is a countable up-final subset of $T$, and there is a countable up-pervasive subset of $V(H)$.
	\end{thm}
	\Proof Let $T^+=\{t\in T:t\ge t_0 \}$, and define $T^-$ similarly. Define $J_{-1}=\emptyset$. Inductively, suppose that $i\ge 0$,
	and $t_0\LL t_i$
	have been defined with the properties that $B_{t_0}\LL B_{t_i}$ are pairwise disjoint, and 
	for all $t\in T$ with $t_0\le t\le t_i$,
	there is at least one and at most two values of $h\in\{0\LL i\}$ with $B_{t_h}\cap B_t\ne\emptyset$.
	Let $J_i$
	be the set of all $t\in T^+$ such that $B_t\cap B_{t_h}\ne \emptyset$ for some $h\in \{0\LL i\}$. Thus $J_i$ is an interval of $T^+$
	containing $t_0$, and $J_0\subseteq J_1\subseteq \cdots\subseteq J_i$. Define $L_i=T^+\setminus J_i$.
	If $L_i=\emptyset$, then the finite set $B_{t_0}\cupcup B_{t_i}$ is up-pervasive, contradicting that $(B_t:t\in T)$ is upper-open.
	So $L_i\ne \emptyset$. Let $Y_{i+1}$ be the set of all southern border
	vertices of $L_i$. By \ref{connected}, $Y_{i+1}\ne \emptyset$,
	and there exists $t_{i+1}\in L_i$ such that $Y_{i+1}\subseteq B_{t_{i+1}}$. Since $t_{i+1}\notin J_i$ it follows that
	$B_{t_{i+1}}\cap B_{t_h}=\emptyset$ for all $h\in \{0\LL i\}$.
	This completes the inductive definition. We see that the sets $B_{t_i}\;(i\ge 0)$ are pairwise disjoint.
	\\
	\\
	(1) {\em For each $t\in T$, there are at most two values of $i\ge 0$ such that $B_t\cap B_{t_i}\ne \emptyset$.}
	\\
	\\
	If $h,j\ge 0$ with $j\ge h+2$, and $t\in T$, we claim that one of $B_t\cap B_{t_h}, B_t\cap B_{t_j}$ is empty.
	To see this, choose $i$ with $h<i<j$. If $t\le t_{i}$, then $B_t\cap B_{t_j}\subseteq B_{t_i}$ from the
	definition of a line-decomposition, and yet $B_{t_{i}}\cap B_{t_j}=\emptyset$, and therefore $B_t\cap B_{t_j}=\emptyset$;
	Similarly if $t_i\le t$ then $B_t\cap B_{t_h}\subseteq B_{t_i}$, and so $B_t\cap B_{t_h}=\emptyset$.
	This proves (1).
	\\
	\\
	(2) {\em $T^+$ equals the union of the sets $J_i\;(i\ge 0)$.}
	\\
	\\
	Let $J$ be the union of the intervals $J_i\;(i\ge 0)$;
	so $J$ is an initial interval of $T^+$. Suppose that $J\ne T^+$.
	Let $L=T^+\setminus J$, and let $Y$ be the set of southern border vertices of $L$. By \ref{connected}, there exists $y\in Y$,
	and there exists $t\in L$ with $Y\subseteq B_t$.
	Choose $s\in J\cap \tau(y)$, and choose $i\ge 0$ with $s\in J_i$. Since $\tau(y)\cap L\ne \emptyset$ and therefore
	$\tau(y)\cap L_{i}\ne \emptyset$, it follows that $y\in Y_{i+1}\subseteq B_{t_{i+1}}$, and so $t\in J_{i+1}$, contradicting that
	$t\notin J$. This proves (2).
	
	\bigskip
	It follows from (1) and (2) that the sequence $t_i\;(i\ge 0)$ satisfies the first two bullets of the theorem. For the third,
	let $t\in T$. From (2) there exists $i\ge 0$ with $t\in J_i$. For each integer $j>i$, since $t_j\notin J_i$ it follows that
	$t<t_j$. This proves the third bullet, and so proves the first assertion of the theorem.
	
	For the remainder, we observe first that since the sequence $\{y_i:i\ge 0\}$ is infinite, the third bullet of the theorem implies that
	$\{y_i:i\ge 0\}$ is up-final.
	Now let $X$ be the union of all the sets $B_{t_i}(i\ge 0)$. Then $X$ is countable, since
	each $B_{t_i}$ is finite; and $X\cap B_t\ne \emptyset$ for each $t\ge t_0$ (by the second bullet of the theorem). To show that $X$
	is up-pervasive, it remains to show that for each $t\in T$, only finitely many vertices of $X$ belong to $\bigcup_{t'\le t}B_{t'}$.
	Let $t\in T$, and choose $j\ge 0$ such that $t_i>t$
	for all $i\ge j$ (this is possible by the third assertion of the theorem).
	Since at most two values of $i$ satisfy $B_{t_i}\cap B_t\ne \emptyset$,
	it follows that $B_i\cap B_t=\emptyset$ for $i\ge j+2$, and therefore
	$B_i\cap \bigcup_{t'\le t}B_{t'}=\emptyset$ for $i\ge j+2$. Consequently, if $x\in X\cap \bigcup_{t'\le t}B_{t'}$, then
	$x$ belongs to one of $B_0\LL B_{j+1}$, and hence the number of such $x$ is finite.
	This proves \ref{getupfinal}.~\bbox
	
	We observe:
	\begin{thm}\label{subset}
		If $(B_t:t\in T)$ is a non-null line-decomposition of a graph $H$, and $X\subseteq V(G)$ is {\em up-pervasive}, then every infinite
		subset of $X$ is up-pervasive.
	\end{thm}
	\Proof
	Let $X'\subseteq X$ be infinite. We need to show that for each $t\in T$, there exists $t'\in T$ with $t'\ge t$ such that
	$X'\cap B_{t'}\ne \emptyset$. Suppose not; and so $t'<t$ for every $t'\in T$ with $X'\cap B_{t'}\ne \emptyset$. In particular,
	$X'\subseteq \bigcup_{t'\le t}B_{t'}$.  But
	since $X$ is up-pervasive, $X\cap \bigcup_{t'\le t}B_{t'}$ is finite, contradicting that $X'$ is infinite. This proves \ref{subset}.~\bbox
	

	The next result uses the up-pervasive and down-pervasive sets  just constructed to find an appropriate geodesic $P$.
	
	\begin{thm}\label{getneargeo}
		Let $(B_t:t\in T)$ be a non-null line-decomposition of a connected graph $H$, with finite width, and let $C\ge 2$.
		Let $S\subseteq T$ be countably infinite.
		Let $\phi$ be a surjective $(C-1,C)$-quasi-isometry from a graph $G$ to $H$. Then for each $s\in S$,
		there exists $v_s\in V(G)$, with $\phi(v_s)\in B_s$, such that for every finite subset $X\subseteq S$, there is a geodesic $P$ in $G$ such that
		$\dist_G(v_s,P)\le C^2-1$ for each $s\in X$.
	\end{thm}
	\Proof Since $\phi$ is surjective, for each $v\in V(H)$ there exists $u\in V(G)$ such that $\phi(u)=v$; choose some such $u$ and denote it by
	$\psi(v)$, for each $v\in V(H)$.
	By \ref{getupfinal}, there is a singleton or countably infinite set $X\subseteq V(H)$ that is up-pervasive.
	Let $Y=\{\psi(x):x\in X\}$.
	It follows that $X,Y$
	have the same cardinality.
	Similarly there is a
	singleton or countably infinite set $Z\subseteq V(G)$ such that $\phi(Z)$ is down-pervasive.
	Since $S$ is countably infinite, we may write $S=\{s_i:i\ge 1\}$.
	
	Take a well-order $\lambda$ of the set of all edges of $G$ (this is possble from the well-ordering theorem). We call $\lambda$ a {\em tie-breaker}.
	If $P,Q$ are distinct paths finite of $G$, we say $P$ is {\em $\lambda$-shorter} than $Q$
	if either
	\begin{itemize}
		\item $|E(P)|<|E(Q)|$; or
		\item $|E(P)|=|E(Q)|$, and the first element (under $\lambda$) of $(E(P)\setminus E(Q))\cup (E(Q)\setminus E(P))$ belongs to $P$.
	\end{itemize}
	This defines a total order on the set of all finite paths of $G$.
	A {\em $\lambda$-geodesic} means a finite path $P$ such that no other path joining its ends is $\lambda$-shorter than $P$.
	Every $\lambda$-geodesic of $G$ is a geodesic of $G$, but the converse is false.
	(The point of the tie-breaker is that there is only
	one $\lambda$-geodesic between any two vertices, while this is not true for geodesics; this will be convenient.)
	It is easy to check that
	if $P$ is a $\lambda$-geodesic then so are all subpaths of $P$.
	Since $H$ is connected and $\phi$ is a quasi-isometry, it follows that $G$ is connected.
	For every two vertices $u,v\in V(G)$, 
	let $P_{u,v}$ be the $\lambda$-geodesic in $G$ between $u,v$.
	
	Let $k\ge 0$, and let $v_i\in V(G)$ for $1\le i\le k$. We say $(v_1\LL v_k)$ is a {\em good choice} if
	\begin{itemize}
		\item $\phi(v_i)\in B_{s_i}$ for $1\le i\le k$; 
		\item there is a set $Y_k\subseteq Y$, either countably infinite or a singleton, such that $\{\phi(y):y\in Y_k\}$ is up-pervasive; and there
		is a set $Z_k\subseteq Z$, either countably infinite or a singleton, such that
		$\{\phi(z):z\in Z_k\}$ is down-pervasive; 
		\item $\dist_G(v_i, V(P_{y,z}))\le C^2-1$ for all $i\in \{1\LL k\}$, all $y\in Y_k$ and all $z\in Z_k$.
	\end{itemize}
	\noindent
	(1) {\em If $k\ge 0$ and $(v_1\LL v_k)$ is a good choice, then then there exists $v_{k+1}$ such that $(v_1\LL v_k, v_{k+1})$ is a good choice.}
	\\
	\\
	If $Y_k$ is infinite,
	there are only
	finitely many vertices $y\in Y_k$ such that $\phi(y)\in \bigcup_{t\le s_{k+1}}B_{t}$, and we may remove them from $Y_k$ by \ref{subset};
	so we may assume that
	$\phi(y)\notin \bigcup_{t\le s_{k+1}}B_{t}$, for each $y\in Y_k$. Thus (even if $Y_k$ is not infinite and hence is a singleton),
	for each $y\in Y_k$,
	there exists $t''\ge s_{k+1}$
	such that $\phi(y)\in B_{t''}$.
	Similarly we may assume that for each $z\in Z_k$, there exists $t'\le s_{k+1}$
	such that $\phi(z)\in B_{t'}$.

	Let $W$ be the set $\{\psi(v):v\in B_{s_{k+1}}\}$.
	Let $y\in Y_k$ and $z\in Z_k$, and choose $t'\le s_{k+1}\le t''$ such that $\phi(y)\in B_{t''}$ and $\phi(z)\in B_{t'}$.
	By \ref{shortjump}, there exists $u\in V(P_{y,z})$ and $v\in B_{s_{k+1}}$ such that $\dist_H(\phi(u),v)\le C-1$, and so
	$\dist_G(u, \psi(v))\le (C-1)^2+C\le C^2-1$.
	Define $\omega(y,z)=\psi(v)$. Thus, for all $y\in Y_k$ and $z\in Z_k$, we have defined $\omega(y,z)\in W$, and
	$\dist_G(P_{y,z},\omega(y,z))\le C^2-1$.
	
	We claim that there exist $v_{k+1}\in W$ and $Y_{k+1}\subseteq Y_k$ and $Z_{k+1}\subseteq Z_k$ such that
	\begin{itemize}
		\item $\omega(y,z) = v_{k+1}$
		for all $y\in Y_{k+1}$ and $z\in Z_{k+1}$; and
		\item $\phi(Y_{k+1})$ is up-pervasive and $\phi(Z_{k+1})$ is down-pervasive.
	\end{itemize}
	To see this, there are four cases. We recall that $|W|$ is finite.
	If $|Y_k|=|Z_k|=1$, let $Y_{k+1}=Y_k=\{y\}$ and $Z_{k+1}=Z_k=\{z\}$ and let $v_{k+1}=\omega(y,z)$.
	If $Y_k=\{y\}$ and $Z_k$ is infinite, then there is an infinite subset $Z_{k+1}$ of $Z_k$
	such that for all $z\in Z_{k+1}$ the vertices $\omega(y,z)$ are all equal (to some $v_{k+1}$); and
	$\{\phi(z)\;(z\in Z_{k+1})\}$ is down-pervasive by \ref{subset}. Similarly the result holds if $Y_k$ is infinite and $|Z_k|=1$. Finally, if both
	$Y_k, Z_k$ are infinite, by an infinite form of Ramsey's theorem for bipartite graphs, there are infinite subsets $Y_{k+1}\subseteq Y_k$
	and $Z_{k+1}\subseteq Z_k$ satisfying the first bullet for some choice of $v_{k+1}$; and again $\{\phi(y)\;(y\in Z_{k+1})\}$ is up-pervasive
	and $\{\phi(z)\;(z\in Z_{k+1})\}$ is down-pervasive, by \ref{subset}. This proves (1).

	\bigskip
	
	Let $v_1,v_2\LL $ be the sequence given by (1); this proves \ref{getneargeo}.~\bbox

	
	\begin{thm}\label{tidyneargeo}
		Let $(B_t:t\in T)$ be a non-null line-decomposition of a graph $H$, with finite width, and let $C\ge 1$.
		Let $\phi$ be a surjective $(C-1,C)$-quasi-isometry from a graph $G$ to $H$.
		Then there is an integer interval $I$, and $v_i\in V(G)$ for each $i\in I$, and an integer $d_{i}>0$
		for all $i\in I$ with $i+1\in I$,
		with the following properties, where $R$ denotes the set of vertices $v$ of $G$ such that $\dist_G(v,\{v_i:i\in I\})\le 3C^2$,
		and $\phi(R)$ denotes $\{\phi(v):v\in R\}$:
		\begin{itemize}
			\item for each $t\in T$ there exist $t', t''\in T$ with $t'\le t\le t''$ such that $\phi(R)\cap B_{t'}\ne \emptyset\ne \phi(R)\cap B_{t''}$;
			\item for all distinct $h,j\in I$, $\dist_G(v_h,v_{j})\ge 2C^2$; and
			\item  for all distinct $h,j\in I$ with $j>h$, $|\dist_G(v_h,v_j)- \sum_{h\le i<j}d_i|\le 2C^2$.
		\end{itemize}
	\end{thm}
	
	\Proof
	If $(B_t:t\in T)$ is not upper-open, choose $y_0\in V(G)$ such that $\{\phi(y_0)\}$ is upfinal, and if
	$(B_t:t\in T)$ is not lower-open, choose $z_0\in V(G)$ such that $\{\phi(z_0)\}$ is downfinal.
	Suppose first that
	$y_0, z_0$ both exist and $\dist_G(y_0,z_0)<2C^2$. Then we may set $I=\{1\}$ and $v_1 = y_0$ and the theorem is satisfied.
	Next, if $y_0,z_0$ both exist and $\dist_G(y_0,z_0)\ge 2C^2$, set $I=\{1,2\}$, $v_1 = y_0$, $v_2 = z_0$ and $d_1 = \dist_G(y_0,z_0)$;
	then again the theorem is satisfied. So we may assume that not both $y_0, z_0$ exist, and from the symmetry we may assume that $y_0$ does not exist, and so $(B_t:t\in T)$ is upper-open.
	
	By \ref{getupfinal}, there is a countable subset of $T$ that is upfinal, and similarly one that is downfinal, either infinite or the
	singleton $\{z_0\}$.
	Let $S_1$
	be their union. Thus $S_1$ is a countable subset of $T$, and satisfies:
	\begin{itemize}
		\item for each $t\in T$, there exist $s,s'\in S_1$ with $s\le t\le s'$; and
		\item for all $t,t'\in T$ with $t<t'$,  there are only finitely many $s\in S_1$ with $t<s<t'$.
	\end{itemize}
	By \ref{getneargeo}, for each $s\in S_1$,
	there exists $v_s\in V(G)$, with $\phi(v_s)\in B_s$, such that for every finite subset $X\subseteq S_1$, there is a geodesic $P$ in $G$ such that
	$\dist_G(v_s,P)\le C^2-1$ for each $s\in X$.
	Choose $S_2\subseteq S_1$,  with $z_0\in S_2$ if $z_0$ exists, maximal such that $\dist_G(v_s,v_{s'})\ge 2C^2$ for all distinct $s,s'\in S_2$ (this is possible by Zorn's lemma).
	
	The next claim is aimed towards the first bullet of the theorem.
	\\
	\\
	(1) {\em  For each $t\in T$, either there exists $s\in S_2$ such that $\dist_G(v_s,\phi^{-1}(B_t))\le 3C^2-1$,
		or there exist $s,s'\in S_2$ with $s<t<s'$.}
	\\
	\\
	Suppose that there is no $s\in S_2$ with $s\ge t$ (the proof is similar if there is no $s$ with $s\le t$).
	From the definition of $S_1$, there exists
	$s'\in S_1$ such that $t\le s'$. Thus $s'\notin S_2$,
	and so there exists $s\in S_2$ such that $\dist_G(v_{s},v_{s'})< 2C^2$ (this exists from the maximality of $S_2$).
	By our assumption, $s<t$, and we may assume that $v_s\notin \phi^{-1}(B_t)$, that is, $\phi(v_s)\notin B_t$.
	Let $P$
	be a geodesic of $G$ between $v_s, v_{s'}$. Thus $P$ has length at most $2C^2-1$. If $\dist_G(P, \phi^{-1}(B_t))\le C^2$,
	then $\dist_G(v_s, \phi^{-1}(B_t))\le 3C^2-1$, as required, so we assume that $\dist_G(P, \phi^{-1}(B_t))> C^2$. In particular,
	$\phi(v_{s'})\notin B_t$. Since $\phi(v_s)\in B_s$ and $\phi(v_{s'})\in B_{s'}$ and $s\le t\le s'$,
	and $\phi(v_s),\phi(v_{s'})\notin B_t$, it follows that $\phi(v_s), \phi(v_{s'})$ are in different components of $H\setminus B_t$.
	Choose a minimal subpath $Q$ of $P$
	with ends $v_s,q$ say, such that $\phi(v_s), \phi(q)$ are not in the same component of $H\setminus B_t$. Let $p$ be the neighbour of
	$q$ in $Q$. It follows that $\phi(p)$ does not belong to the component of $H\setminus B_t$ that contains $\phi(q)$. Since
	$\dist_H(\phi(p),\phi(q))\le 2C-1$, it follows that $\dist_H(\phi(P),B_t)\le C-1$, and so $\dist_G(P, \phi^{-1}(B_t))\le C^2-1$.
	Since $P$ has length at most $2C^2$, it follows that $\dist_G(v_s,\phi^{-1}(B_t))\le 3C^2-1$. This proves (1).
	
	\bigskip
	
	Let $X\subseteq S_2$ be finite; then
	there is a geodesic $P$ in $G$ such that
	$\dist_G(v_s,P)\le C^2-1$ for each $s\in X$. For each $s\in X$, choose $p_s\in V(P)$ such that $\dist_G(v_s,p_s)\le C^2-1$. For all
	$s,s'\in X$, let $n(s,s')=\dist_G(p_s,p_{s'})$. We call
	$(n(s,s'):s,s'\in X)$ a {\em gap matrix} for $X$. One set $X$ might have several gap matrices, because
	the matrix also depends on  $P$ and the choices of the vertices $p_s\;(s\in X)$. Since
	$$|\dist_G(p_s,p_{s'})-\dist_G(v_s,v_{s'})|\le 2(C^2-1),$$
	there are fewer than $4C^2$ possibiilities for each entry $n_{s,s'}$ of the gap matrix, and so at most $(4C^2)^{|X|^2}$
	possibilites for the gap matrix for $X$.
	
	If $X$ is a finite subset of $S_2$ and $(n(s,t):s,t\in X)$ is a gap matrix for $X$, let $X'\subseteq X$; then
	$(n(s,t):s,t\in X')$ is a gap matrix for $X'$, and we say $(n(s,t):s,t\in X)$ {\em extends} $(n(s,t):s,t\in X')$.
	Take a sequence $X_1\subseteq X_2\subseteq \cdots$ of finite subsets of $S_2$ with union $S_2$, and for each $i\ge 1$,
	let $N_i$ be a corresponding gap matrix for $X_i$. Make a graph $K$ with vertex set the set of all gap matrices for each $X_i$, in which
	a gap matrix for $X_i$ is adjacent to a gap matrix for $X_j$ if $j=i+1$ and the second gap matrix extends the first. Then $K$ is a rooted 
	tree
	with infinitely many vertices and all degrees finite, and so it has an infinite path, by K\"onig's lemma. Consequently there exists
	$d(s,s')$ for all $s,s'\in S_2$, such that for every finite $X\subseteq S_2$, $(d(s,s'):s,s'\in X)$ is a gap matrix for $X$.

	Let $X=\{x_1,x_2\}\subseteq S_2$ with $|X|=2$; then $(d(s,s'):s,s'\in X)$ is a gap matrix for $X$. Let $P$ and $p_s\;(s\in X)$
	be as in the definition of gap matrix. The vertices $p_{x_1}, p_{x_2}$ are distinct, since the vertices $v_{x_1}, v_{x_2}$
	have distance at least $2C^2>2 (C^2-1)$. Consequently $d(x_1,x_2)\ge 1$ for all distinct $x_1,x_2\in S_2$.
	
	Let $X\subseteq  S_2$ with $|X|=3$; then $(d(s,s'):s,s'\in X)$ is a gap matrix for $X$. Let $P$ and $p_s\;(s\in X)$
	be as in the definition of gap matrix.
	Since the vertices $p_{s}(s\in X)$ all belong to the geodesic $P$,
	one of them is between the other two in $P$; let $X=\{x_1,x_2,x_3\}$ where $p_{x_2}$ is between the other two in $P$.
	It follows that $d(x_1,x_2)+ d(x_2,x_3) = d(x_1,x_3)$. Consequently $d(x_1,x_2), d(x_2,x_3)<d(x_1,x_3)$ since they are all
	strictly positive.
	
	Since this holds for all triples of distinct vertices in $S_2$, there is a linear ordering $<$ of $S_2$ such that if $a<b<c\in X$
	then $d(a,b)+ d(b,c) = d(a,c)$. From the additivity of the $d$ function, there are only finitely many terms between any two terms
	of this linear order, and so we can number $S_2$ as $\{v_i:i\in I\}$, where $I$ is an interval of integers, and
	define $d_i = d(v_i,v_{i+1})$ for each $i\in I$ with $i+1\in I$, such that
	$d_{h,j} = \sum_{h\le i<j}d_i$ for all $h,j\in I$ with $j>h$, and consequently,
	$$|\dist_G(v_h,v_j)- \sum_{h\le i<j}d_i|\le 2C^2.$$
	This proves \ref{tidyneargeo}.~\bbox

	Now we can prove \ref{addgeo0}, which we restate:
	\begin{thm}\label{addgeo}
		Let $(B_t:t\in T)$ be a non-null line-decomposition of a graph $H$, with finite width, and let $C\ge 2$.
		Let $\phi$ be a surjective $(C-1,C)$-quasi-isometry from a graph $G$ to $H$.
		Then there is a graph $G'$ and a geodesic $P$ of $G'$ with the following properties:
		\begin{itemize}
			\item $G$ is an induced subgraph of $G'$, and $\dist_G(u,v) = \dist_{G'}(u,v)$ for all $u,v\in V(G)$;
			\item the identity map from $V(G)$ into $V(G')$ is a $(1,2C^2+1)$-quasi-isometry, and there is a
			$(1,4C^2+2)$-quasi-isometry from $G'$ to $G$ that maps each vertex of $G$ to itself;
			\item for each $t\in T$, $\dist_{G'}(P,\phi^{-1}(B_t))\le 6C^2$.
		\end{itemize}
	\end{thm}
	\Proof By \ref{tidyneargeo}, there is an integer interval $I$, and $v_i\in V(G)$ for each $i\in I$, and an integer $d_{i}>0$
	for all $i\in I$ with $i+1\in I$,
	with the following properties, where $R$ denotes the set of vertices $v$ of $G$ such that $\dist_G(v,\{v_i:i\in I\})\le 3C^2$,
	and $\phi(R)$ denotes $\{\phi(v):v\in R\}$:
	\begin{itemize}
		\item for each $t\in T$ there exist $t', t''\in T$ with $t'\le t\le t''$ such that $\phi(R)\cap B_{t'}\ne \emptyset\ne \phi(R)\cap B_{t''}$;
		\item  for all distinct $h,j\in I$ with $j>h$, $|\dist_G(v_h,v_j)- \sum_{h\le i<j}d_i|\le 2C^2$.
	\end{itemize}
	For each $i\in I$, let $u_i$ be a new vertex.
	For each $i\in I$ with $i+1\in I$, let $Q_i$ be a geodesic of $G$ between $v_i, v_{i+1}$, and let $P_i$ be a path of new vertices with ends $u_i,u_{i+1}$, of length $d_i$. For each such $i$, the lengths of $P_i$ and $Q_i$ differ by at most $2C^2$. Let $P$ be the
	union of the paths $P_i$, over all $i\in I$ with $i+1\in I$. Thus $P$ is a path. Define a map $\alpha$ from $V(P)$ into $V(G)$
	as follows. Let $p\in V(P)$. If $p=u_i$ for some $i\in I$ then $\alpha(p) = v_i$. Otherwise $p$ is an internal vertex of $P_i$
	for some $i\in I$ with $i+1\in I$. Let $P_i[u_i,p]$ have length $h$. If $Q_i$ has length at least $h$,
	let $\alpha(p)$ be the vertex $q\in V(Q_i)$ such that $Q_i[v_i,q]$ has length $h$, and otherwise let $\alpha(p) = v_{i+1}$. This defines $\alpha$.
	
	Now for each $p\in V(P)$, add a path $R_p$ of new vertices between $p, \alpha(p)$ of length $2C^2+1$. Let $G'$ be the union of $G,P$
	and all the paths $R_p\;(p\in V(P))$. We will show that $G'$ and $P$ satisfy the theorem.
	\\
	\\
	(1) {\em For all $p,p'\in V(P)$, let $q=\alpha(p)$ and $q'=\alpha(p')$; then $\dist_P(p,p')$ and $\dist_G(q,q')$ differ by at most $4C^2$.}
	\\
	\\
	Suppose first that
	there exists $i\in I$ such that $p,p'\in V(P_i)$. Then $\alpha(p),\alpha(p')\in V(Q_i)$. If $q,q'\ne v_{i+1}$, then
	$\dist_P(p,p')$ and $\dist_G(q,q')$ are equal, from the definition of $\alpha$. Hence we assume that $q'=v_{i+1}$. If also $q=v_{i+1}$,
	then $P[u_i,p]$ and $P[u_i, p']$ both have length at least $i|E(Q_i)|$, and so
	$\dist_P(p,p')\le |E(P_i)|-|E(Q_i)|\le 2C^2$ and $\dist_G(q,q')=0$, as required. Hence we may assume that $q\ne v_{i+1}$, and therefore
	$P_i[u_i,p]$ has the same length as $Q_i[v_i, q]$. Since $\dist_P(u_i,p')\ge |E(Q_i)|$, it follows that
	$\dist_P(p,p')\ge \dist_G(q,q')$. But
	$$\dist_P(u_i,p')\le |E(P_i)|\le |E(Q_i)|+2C^2.$$
	Since $\dist_P(u_i,p')=\dist_P(u_i,p) + \dist_P(p,p')$, and $\dist_G(v_i,q')=\dist_G(v_i,q)+\dist_G(q,q')$,
	we deduce that $\dist_P(p,p')\le \dist_G(q,q')+4C^2$, as required.
	
	So we may assume that there is no such $i\in I$, and therefore there are vertices of the form $u_i\;(i\in I)$ that are internal vertices
	of $P[p,p']$. We may assume that $p\in V(P_h)\setminus \{u_{h+1}\}$, and $p\in V(P_j)\setminus \{u_j\}$, where $h,j\in I$ with $h<j$,
	and $u_{h+1}\LL u_j$ are all
	internal vertices of $P[p,p']$.
	
	We claim first that $\dist_P(p,p')\le \dist_G(q,q')+4C^2$.
	To see this, observe that 
	$$\dist_G(v_h,q) + \dist_G(q,q')+\dist_G(q', v_{j+1})\ge \dist_G(v_h, v_{j+1})\ge d_h+d_{h+1}+\cdots+d_j-2C^2.$$
	But the latter equals
	$$ \dist_P(u_h,p) + \dist_P(p,p')+\dist_P(p', u_{j+1})-2C^2.$$
	Consequently
	$$\dist_G(v_h,q) + \dist_G(q,q')+\dist_G(q', v_{j+1})\ge \dist_P(u_h,p) + \dist_P(p,p')+\dist_P(p',u_{j+1})-2C^2,$$
	that is,
	$$\dist_G(q,q')-\dist_P(p,p') \ge \dist_P(u_h,p) -\dist_G(v_h,q) + \dist_P(p',u_{j+1})-\dist_G(q', v_{j+1})-2C^2.$$
	But $\dist_P(u_h,p) \ge \dist_G(v_h,q)$, and $\dist_P(p',u_{j+1})-\dist_G(q', v_{j+1})\ge -2C^2$;
	so $\dist_G(q,q')-\dist_P(p,p')\ge -4C^2$, as claimed
	
	Next, we claim that $\dist_P(p,p')\ge \dist_G(q,q')-4C^2$. To see this, observe that 
	$$\dist_G(q,q')\le \dist_G(q,v_{h+1}) + \dist_G(v_{h+1}, v_j)+ \dist_G(v_j,q').$$
	But $\dist_G(q,v_{h+1})\le \dist_P(p,u_{h+1})+2C^2$, and
	$$\dist_G(v_{h+1}, v_j)\le d_{h+1}+\cdots d_{j-1} +2C^2=\dist_P(u_{h+1}, u_j)+2C^2,$$
	and
	$\dist_G(v_j,q')\le \dist_P(u_j,p')$. Consequently
	$$\dist_G(q,q')\le \dist_P(p,u_{h+1})+2C^2 + \dist_P(u_{h+1}, u_j)+2C^2+ \dist_P(u_j,p')=\dist_P(p,p')+4C^2,$$
	as claimed. 
	This proves (1).
	\\
	\\
	(2) {\em No finite geodesic of $G'$ with both ends in $V(G)$ has a vertex in $P$. Consequently
		if $u,v\in V(G)$, then $\dist_G(u,v) = \dist_{G'}(u,v)$.}
	\\
	\\
	The second statement follows immediately from the first. Suppose that the first is false, and let $L$ be the shortest geodesic
	of $G'$ with both ends in $V(G)$ and with a vertex in $V(P)$. It follows that the ends of $L$ are $\alpha(p), \alpha(p')$
	for some $p,p'\in V(P)$; and $L$ is the union of $R_p, R_{p'}$ and the subpath of $P$ between $p,p'$.  Thus $L$ has length
	$$2(2C^2+1) +\dist_P(p,p')\ge 2+\dist_G(\alpha(p), \alpha(p'))\ge 2+|E(L)|$$
	since $L$ is a geodesic of $G'$, a contradiction. This proves (2).
	
	\bigskip
	Since every vertex in $V(G')$ has distance in $G'$ at most $2C^2+1$ from some vertex of $G$, (2) implies that the identity map is a
	$(1,2C^2+1)$-quasi-isometry from $G$ to $G'$. By (1), the map $\beta$ from $V(G')$ to $V(G)$, with $\beta(v)=v$ for $v\in V(G)$,
	and with $\beta(v)=\alpha(p)$ for each $p\in V(P)$ and $v\in V(R_p)$, is a $(1,4C^2+2)$-quasi-isometry from $G'$ to $G$.
	
	Finally, let $t\in T$. From the application of \ref{tidyneargeo},
	there exist $t', t''\in T$ with $t'\le t\le t''$ such that $\phi(R)\cap B_{t'}\ne \emptyset\ne \phi(R)\cap B_{t''}$;
	where $R$ denotes the set of vertices $v$ of $G$ such that $\dist_G(v,\{v_i:i\in I\})\le 3C^2$,
	and $\phi(R)$ denotes $\{\phi(v):v\in R\}$. Since $\phi(R)\cap B_{t'}\ne \emptyset$, there exists $h\in I$ and $x_h\in V(G)$
	such that $\dist_G(v_h, x_h)\le 3C^2$ and $\phi(x_h)\in B_{t'}$.  Let $F_h$ be a geodesic in $G$ between $v_h, x_h$.
	Similarly  there exists $j\in I$ and $x_j\in V(G)$
	such that $\dist_G(v_j, x_j)\le 3C^2$ and $\phi(x_j)\in B_{t''}$, and a geodesic $F_j$ between $v_j, x_j$.
	We may assume that $h\le j$.
	The union of the paths $F_h, Q_h, Q_{h+1}\LL Q_{j-1}, F_j$
	is a connected subgraph of $G$ containing $x_h$ and $x_j$, and so one of its vertices has distance in $G$ at most $C^2-1$ from
	$\phi^{-1}(B_t)$, by \ref{shortjump}. But each of its vertices has distance at most $3C^2+(2C^2+1)$ from $P$ in $G'$, since the vertices of $F_h$
	have distance at most $3C^2$ from $v_h$ and hence at most $3C^2+(2C^2+1)$ from $u_h$, and the same holds for $F_j$, and for
	$h\le i\le j-1$, each vertex of $Q_i$ has distance in $G'$ at most $2C^2+1$ from $P$.  Hence $\dist_{G'}(P,\phi^{-1}(B_t))\le 6C^2$.
	This proves \ref{addgeo}.~\bbox


\end{document}